\documentclass{amsart}
\usepackage[T1]{fontenc}
\usepackage[english]{babel}
\usepackage{amsmath,amsthm,amssymb, latexsym,geometry,stackrel,enumerate}
\usepackage[all]{xy}

\usepackage[dvips]{graphicx}
\DeclareGraphicsExtensions{.bmp,.png,.pdf,.jpgg}
\usepackage[all]{xy}
\usepackage{verbatim}
\usepackage[mathcal]{euscript}
\usepackage{epsfig}
\usepackage{multicol}
\usepackage{color}
\usepackage{verbatim}
\usepackage{graphicx,float}

\swapnumbers
\theoremstyle{plain}
\newtheorem{teo}{Theorem}[section]
\newtheorem{lema}[teo]{Lemma}
\newtheorem{prop}[teo]{Proposition}
\newtheorem{coro}[teo]{Corollary}

\theoremstyle{definition}
\newtheorem{defi}[teo]{Definition}
\newtheorem{obs}[teo]{Remark}

\newtheorem{ejem}[teo]{Example}
\newtheorem{algo}[teo]{Algorithm}

\newcommand{\ts}[1]{\normalfont{\textsf{#1}}}
\newcommand{\K}{\ts k}

\renewcommand{\a}{\alpha}
\renewcommand{\b}{\beta}
\newcommand{\e}{\varepsilon}

\renewcommand{\d}{ \delta}
\newcommand{\s}{\sigma}

\newcommand{\BBeq}{tilting-cotilting equivalent  }

\def\Ho{\mathop{\rm Hom}\nolimits}

\def\En{\mathop{\rm End}\nolimits}

\newcommand{\mor}[3]{$#1\colon #2 \to #3$}

\title[ Derived class of $m$-cluster tilted algebras of type $\widetilde{\mathbb{A}}$]{ Derived class of $m$-cluster tilted algebras of type $\widetilde{\mathbb{A}}$}

\author[V. Gubitosi]{Viviana Gubitosi}
\address{Instituto de Matem\'{a}tica y Estad\'{\i}stica Rafael Laguardia, Facultad de Ingenier\'{\i}a - UdelaR, Montevideo, Uruguay, 11200 }
\email{gubitosi@fing.edu.uy}

\begin{document}
\maketitle

\begin{abstract}
 In this paper, we characterize all the finite dimensional algebras that are derived equivalent to an  $m-$cluster tilted algebras of type $\widetilde{\mathbb{A}}$. This generalizes a result of Bobi\'nski and Buan \cite{BB10}.
\end{abstract}

\section*{Introduction}

Cluster categories were introduced in \cite{Buan2006} as a representation theoretic framework for the  cluster algebras of Fomin and Zelevinski \cite{FZ02}. The clusters correspond to the tilting objects in the cluster category.

Given an hereditary finite dimensional algebra $H$  over  an algebraically closed field $\K$ the $m$-cluster category is defined to be $\mathcal{C}_m(H):=\mathcal{D}^b(H)/ \tau^{-1} [m]$, where $[m]$  denotes the $m$-th power of  the shift functor $[1]$ and $\tau$ is the Auslander - Reiten translation in $\mathcal{D}^b(H)$. By a result of Keller \cite{K05}, the $m$-cluster category is triangulated.  For the  $m$-cluster category,  $m$-cluster tilting objects have been defined by Thomas, in \cite{Thomas2007}, who in addition showed that they are in bijective correspondence with the $m$-clusters associated by  Fomin and Reading to a finite root system in \cite{FR05}. The endomorphism algebras of the $m$-cluster tilting objects are called $m$-cluster tilted algebras or, in case $m=1$, cluster tilted algebras.

In \cite{ABCP09}, Assem \emph{et al.} showed that cluster tilted algebras coming from triangulations of the disc or the annulus with marked points on their boundaries are gentle, and, in fact, that these are the only gentle cluster tilted algebras. The class of gentle algebras  defined by Assem and Skowro\'nski in \cite{AS87}  has been extensively studied, see \cite{AH81, Avella-Alaminos2008, BB10, Buan2008, Murphy2010, SZ03}, for instance, and is particularly well understood, at least from the representation theoretic point of view. This class includes, among others, iterated tilted and cluster
tilted algebras of types $\mathbb{A}$ and $\tilde{\mathbb{A}}$, and, as shown in \cite{SZ03}, is closed under derived equivalence.

When studying module categories, one is often interested in them up to derived equivalence, or tilting-cotilting equivalence. In \cite{Buan2008}, Buan and Vatne gave a criterion to decide whether two cluster tilted algebras of type $\mathbb{A}$ are themselves derived equivalent or not. This has been done using the determinant of the Cartan matrix as derived invariant, as well as mutations of quivers. Later, Bastian, in \cite{Bastian09}, gave an analogous classification for the $\tilde{\mathbb{A}}$ case. She used another thinner derived invariant, the function $\phi$ introduced by Avella-Alaminos and Geiss in \cite{Avella-Alaminos2008}. In \cite{BB10} a more general question has been considered, namely the characterization of the algebras that are derived equivalent to cluster tilted algebras of type $\mathbb{A}$ or $\tilde{\mathbb{A}}$. Again, in this paper the map $\phi$ is of central importance, and the characterizations therein are given in terms  of this map. In  another direction, results analogous to those of \cite{Buan2008} have been established for $m$-cluster tilted algebras of type $\mathbb{A}$ by Murphy in \cite{Murphy2010}: he described these algebras by quivers and relations, and gave a criterion permitting to decide whether two $m$-cluster tilted algebras of type $\mathbb{A}$ are derived equivalent or not. Again, he used the determinant of the Cartan matrix as in \cite{Buan2008}, but ``elementary polygonal moves'' instead of  -- but equivalently to -- mutations. Later, Bustamante and Gubitosi, in \cite{GubBust}, classified the algebras that are derived equivalent to $m$-cluster tilted algebras of type $\mathbb{A}$  using the Hochschild cohomology ring as derived invariant.

The aim of this paper is to classify the algebras that are derived equivalent to $m$-cluster tilted algebras of type $\mathbb{\widetilde{A}}$. Since it will be turn out that the algebras we are interested in are gentle, we can use the function $\phi$ as a derived invariant.

We now state the main result of this paper (for the definitions of the terms used, we refer the reader to sections 1 and 3 below).

\subsection*{Theorem A}\textit{ Let $Q$ be  a quiver having a root cycle.  A connected algebra $A=\K Q/I$ is derived equivalent to a connected component of an $m$-cluster tilted algebra of type $\mathbb{\widetilde{A}}$ if and only if  $A$ is $\tilde{\mathbb{A}}$-branched. }

\medskip

In particular, specializing to the case $m=1$, we recover known results of \cite{BB10}, and we obtain a criterion allowing to decide whether or not an algebra is derived equivalent to a cluster tilted algebra of type $\mathbb{\widetilde{A}}$. The latter is very easy to use, as it does not require any computation, in contrast with the known result of \cite{BB10}. Until now it was not known whether the function $\phi$ is a complete derived  invariant in general, here, we prove that this is not the case by furnishing a counterexample in the case of algebras derived equivalent to an $m$-cluster tilted algebra of type $\mathbb{\widetilde{A}}$.\\


The paper is organized as follows: In section 1 we recall facts about gentle algebras, derived and tilting-cotiling equivalences, and $m$-cluster tilted algebras. Also we recall what Brenner-Butler tilting modules  and the Avella-Alaminos-Geiss map are. In section 2 we establish the facts about $m$-cluster tilted algebras of type $\mathbb{\widetilde{A}}$ that will be used in the sequel.  In section 3 we introduce what we call $\widetilde{\mathbb{A}}$-branched algebras, precisely those named in the theorem above, and in section 4 we introduce what normal forms are (two particular classes of $\widetilde{\mathbb{A}}$-branched algebras).  In section 5 we start the procedure to reduce an  $\widetilde{\mathbb{A}}$-branched algebra to a normal form. Section 6  and 7 are devoted to the proof of the main theorem and some consequences, among which  we recover the known results mentioned above.

\section{Preliminaries}
\subsection{Gentle algebras}

While we briefly recall some  concepts concerning bound quivers and algebras, we refer the reader to \cite{ASS06} or \cite{ARS95}, for instance, for unexplained notions.

Let \K \  be a commutative field. A quiver $Q$ is the data of two sets, $Q_0$ (the \textit{vertices}) and $Q_1$ (the \textit{arrows}) and two maps \mor{s,t}{Q_1}{Q_0} that assign to each arrow $\a$ its \textit{source} $s(\a)$ and its \textit{target} $t(\a)$. We write \mor{\a}{s(\a)}{t(\a)}. If $\b\in Q_1$ is such that $t(\a)=s(\b)$ then the composition of $\a$ and $\b$ is the path $\a\b$. This extends naturally to paths of arbitrary positive length. The \emph{path algebra} $\K Q$ is the $\K$-algebra whose basis is the set of all paths in $Q$, including one stationary path $e_x$ at each vertex $x\in Q_0$, endowed with the  multiplication induced from the composition of paths. In case $|Q_0|$ is finite, the sum of the stationary paths  - one for each vertex - is the identity.

If the quiver $Q$ has no oriented cycles, it is called \emph{acyclic}. A \emph{relation} in $Q$ is a $\K$-linear combination of paths of length at least $2$ sharing source and target.  A relation which is a path is called \emph{monomial}, and the relation is \emph{quadratic} if all the  paths appearing in it have length $2$. Let $\mathcal{R}$ be a set of relations.
 Let $\langle Q_1 \rangle$ denote the two-sided ideal of $\K Q$ generated by the arrows, and $I$ be the one generated by $\mathcal{R}$. Then $I\subseteq \langle Q_1\rangle ^2$.
 The ideal $I$ is called \emph{admissible} if there exists a natural number $r\geqslant 2$ such that $\left\langle Q_1 \right\rangle^r \subseteq I$. The pair $(Q,I)$ is a \emph{bound quiver}, and associated to it is the algebra $A=\K Q/I$.
It is known that any finite dimensional basic and connected algebra over an algebraically closed field is obtained in this way, see \cite{ASS06}, for instance.

Recall from \cite{AS87} that an algebra  $A= \K Q/I$ is said to be \emph{gentle} if  $I=\langle \mathcal{R} \rangle $, with $\mathcal{R}$ a set of monomial quadratic relations such that :
\begin{enumerate}
 \item[G1.] For every vertex $x\in Q_0$ the sets $s^{-1}(x)$ and $t^{-1}(x)$ have cardinality at most two;
 \item[G2.] For every arrow $\a\in Q_1$ there exists at most one arrow $\b$ and one arrow $\gamma$ in $ Q_1$ such that $\a\b\not\in I$, $\gamma\a\not\in I $;
 \item[G3.] For every arrow $\a\in Q_1$ there exists at most one arrow $\b$ and one arrow $\gamma$ in $Q_1$ such that $\a\b\in I$, $\gamma\a\in I$.
\end{enumerate}

Gentle algebras are special biserial (see \cite{WW85}), and have extensively been studied in several contexts, see for instance \cite{Avella-Alaminos2008, BB10, Buan2008, Murphy2010, SZ03}.

\subsection{Tilting-cotilting, and derived equivalences}
Given a finite dimensional algebra $A=\K Q/I$ a \emph{tilting} module is a finitely generated right $A$-module of projective dimension less than or equal to 1, having no self extensions and exactly $|Q_0|$ indecomposable non isomorphic direct summands, \cite{ASS06}. The notion of cotilting module is defined dually. Given a tilting $A$-module $T$, with $A$ hereditary, the algebra $\En_A{T}$ is said to be \emph{tilted}. Two algebras $A$ and $B$ are said to be \emph{tilting-cotilting} equivalent if there exists a finite sequence of algebras $A=A_0, A_1,\ldots, A_r=B$ and $A_i$-tilting (or cotilting) modules $T_i$ such that $A_{i+1}=\En_{A_i}{T_i}$ for $i\in\{0,1,\ldots, r-1\}$.

Denote by $\mathcal{D}^b(A)$ the bounded derived category of finite dimensional right $A$-modules. Its objects are bounded complexes of finite dimensional right $A$-modules, and morphisms are obtained from morphisms of complexes by localizing with respect to quasi-isomorphisms (see \cite{H88}). The category $\mathcal{D}^b(A)$ is triangulated, with translation functor induced by the shift of complexes. Two algebras $A$ and $B$ are \emph{derived equivalent} if the categories $\mathcal{D}^b(A)$ and $\mathcal{D}^b(B)$ are equivalent as triangulated categories. It has been shown by Happel \cite{H88} that if two algebras are tilting-cotilting equivalent, then they are derived equivalent. Moreover, Schr\"oer and Zimmermann showed in \cite{SZ03} that the class of gentle algebras is stable under derived equivalence.

\subsection{Brenner  - Butler tilting modules}\label{mutaciones}
Let $(Q,I)$ be a gentle bound quiver without loops and $x\in Q_0$ such that whenever there is an arrow $\alpha$ leaving $x$,   then there is an arrow $\beta$ entering $x$, such that $ \a\b\not\in I$. This includes for instance the vertices that are not the source of any arrow, but excludes the sources of $Q$. Since $(Q,I)$ is gentle, the vertex $x$ has at most two arrows leaving it, say $\a_0$ and $\a_1$. If this is the case, let $\b_0, \b_1$ be the arrows such that $\a_i\b_i \not\in I$, for $i\in\{0,1\}$. Read indices modulo $2$, then for each $i$, there exists at most one arrow $\gamma_{i+1}$ such that $\b_i\gamma_{i+1}  \in I$. Note that since the algebra is gentle, we have $\b_{i+1}\a_i ,\ \b_i\gamma_{i+1} \in I$ (see the left figure on next page).

In \cite{AH81} Assem and Happel showed that an algebra whose quiver is a gentle tree is tilting-cotilting equivalent to an hereditary algebra of type $\mathbb{A}$. This had been done by explicitly giving a sequence of tilting and cotilting modules. At each stage the gentle tree is transformed until the quiver of type $\mathbb{A}$ is reached. We will exhibit an analogous process, called \emph{``elementary transformation over a vertex''} in \cite[Section 7]{AA07}, see also \cite[Section 2]{BB10} or \cite[Section 6]{GubBust}.

\begin{defi}
Let $(Q,I)$ be a gentle bound quiver, and $x$ as above. With these notations the bound quiver obtained by \emph{mutating} $(Q,I)$ at $x$ is the bound quiver defined by $(Q',I')=\sigma_x(Q,I)$ where:
\begin{itemize}
 \item $Q'_0 =Q_0$,
 \item $Q'_1 = Q_ 1 \backslash\{\a_i,\b_i,\gamma_i|i=0,1\}\cup \{\a'_i,\b'_i,\gamma'_i|i=0,1\}$ such that \mor{\a'_i}{b_i}{a_i}, \mor{\b'_i}{x}{b_i}, \mor{\gamma'_i}{c_i}{x}.
\end{itemize}
Let $\mathcal{R}$ be a minimal set of relations generating $I$, in particular it contains $\a_i \b_{i+1},\ \b_i\gamma_{i+1}$. Let $\mathcal{R}'$ be  obtained by replacing in $\mathcal{R}$ the latter by  $\a'_i\b'_{i+1}$, $\b'_i\gamma'_{i+1}$ for $i=0,1$, and, again, indices are to be  read modulo 2. Then $I'$ is the ideal generated by $\mathcal{R}'$.
\end{defi}

 In the sequel, a dotted line joining two arrows 
means, as usual, that their composition belongs to $I$.

\begin{center}
 \begin{tabular}{ccccccc}
 $\SelectTips{eu}{10}\xymatrix@C=.3pc@R=.2pc{    a_0  & && & & & b_1  \ar[dddlll]^{\b_1}  & \ar@/^/@{.}[lld]&& c_0  \ar[lll]_{\gamma_0} \\
&&&&&&&&&\\
&& \ar@/^/@{.}[rr]&&&&&&& \\
&&& x
 \ar[uuulll]^{\a_0} \ar[dddlll]_{\a_1} &&&&&&\\
&&\ar@/_/@{.}[rr]&&&&&&& \\
&&&&&&&&&\\
  a_1  & && & & & b_1  \ar[uuulll]_{\b_0} &\ar@/_/@{.}[llu] && c_1  \ar[lll]^{\gamma_1} }$&  &&& & & $\SelectTips{eu}{10}\xymatrix@C=.3pc@R=.2pc{  a_0 && \ar@/_/@{.}[rrd]& b_1   \ar[lll]_{\a_0'} & && & & & c_0  \ar[dddlll]^{\gamma_0'} \\
&&&&&&&&&\\
&&&&& \ar@/^/@{.}[rr]&&&& \\
&&&&&& x \ar[uuulll]^{\b_1'} \ar[dddlll]_{\b_0'} &&&\\
&&&&&\ar@/_/@{.}[rr]&&&& \\
&&&&&&&&\\
 a_1 & &\ar@/^/@{.}[rru] & b_1  \ar[lll]^{\a_1'} & && & & & c_1  \ar[uuulll]_{\gamma_1'}}$\\
 &  & &&& &\\
$(Q,I)$ & &&& & & $\s_x(Q,I)$.
   \end{tabular}
   \end{center}
\medskip


We have the known result:

\begin{lema}[\cite{BB80}] Let $A=\K Q/I$ be a gentle algebra and $x\in Q_0$ as above. Then
\begin{enumerate}
 \item[$a)$] The module $T_x = \tau^{-1} S_x \oplus \overline{T}=\underset{y\not=x}{\bigoplus} P_y$ is a tilting $A$-module;
 \item[$b)$] The quiver of $\En_{A}(T_x)$ is precisely $\sigma_x(Q,I)$.
\end{enumerate}
\end{lema}

The tilting module $T_x$ is called the \emph{Brenner - Butler tilting module at $x$}, or BB tilting module, for short. In an analogous way one can define the \emph{BB cotilting module at a vertex $y$}, and the corresponding mutation $\s'_y$ on the bound quivers.

\subsection{The Avella Alaminos - Geiss map $\phi$}\label{subsec:AG-phi}
As mentioned  before, in \cite{BB10} the map $\phi$, which is a derived invariant for gentle algebras \cite{Avella-Alaminos2008}, is the main tool used to establish the derived equivalence classification therein. It is a map \mor{\phi}{\mathbb{N} \times \mathbb{N} }{\mathbb{N}} that counts special sequences of paths and relations  in a gentle quiver  $(Q,I)$. We are interested in characterizing this map for normal forms $\widetilde{N}_{n_1,k_1,n_2,k_2,r}$ and $B_{k,n,t}$ defined in section \ref{normal forms}. In what follows we closely follow the exposition of \cite{Avella-Alaminos2008}.

Let $A=\K Q/I$ be a gentle algebra. A \emph{permitted thread} of $A$ is a path $w=\a_1 \a_2\cdots \a_n$ not belonging to $I$, and of maximal length for this property. A \emph{forbidden thread} is a sequence $\pi=\a_n \a_2\cdots \a_1$ formed by pairwise different  arrows with maximal length and such that $\a_{i+1}\a_i\in I$ for all $i\in \{1,2,\cdots, n-1\}$.

We also need trivial permitted and forbidden threads. Let $x\in Q_0$ be such that the sets $s^{-1}(x)$ and $t^{-1}(x)$ have both cardinality at most one. The stationary path at $x$ is a \emph{trivial permitted thread} if when $\alpha$ ending at $x$ and  $\beta$ starting in $x$ are two arrows, then $\b\a\not\in I$. We denote this thread by $h_x$. Similarly, the stationary path at $x$ is a \emph{trivial forbidden thread} if when $\alpha$ ending at $x$ and  $\beta$ starting in $x$ are two arrows, then $\b\a\in I$. We denote by $p_x$ this thread. Assume $x$  and $y$ are two vertices such there is only one  arrow $\a$ entering $x$, an arrow \mor{\b}{x}{y}, and only one arrow $\gamma$ leaving $y$. If both $\b\a$, $\gamma\b\in \mathcal{R}$, then $\b$ is a permitted thread, whereas in case $\b\a$, $\gamma\b\not\in \mathcal{R}$ the arrow $\b$ is a forbidden thread.

Given that $(Q,I)$ is gentle, from \cite{BR87} one knows that there exist maps \mor{\s , \e}{Q_1}{\{\pm 1\}} satisfying:
\begin{itemize}
 \item $\s(\b_0) = -\s(\b_1)$ whenever $\b_0$ and $\b_1$ are arrows sharing their source;
 \item $\e(\b_0) = -\e(\b_1)$ whenever $\b_0$ and $\b_1$ are arrows sharing their target;
 \item If $\b\a $ is a path not belonging to $I$, then $\s(\b)= - \e(\a)$.
\end{itemize}
These maps, that one can set \emph{``quite arbitrarily''}, as noted in \cite[p. 158]{BR87},  extend to paths, thus to threads: given $w=  \a_n\cdots\a_2\a_1$, set $\s(w)=\s(\a_1)$ and $\e(w)=\e(\a_n)$. We extend this to trivial threads as follows: If $h_x$ is a trivial permitted thread, the connectedness of $Q$ assures the existence of an arrow  $\gamma$ leaving $x$ or an arrow  $\b$ ending at $x$. Then in the first case  put $\s(h_x)=-\e(h_x)=-\s(\gamma)$, for the second case put $\s(h_x)=-\e(h_x)=\e(\b)$. Similarly, if $p_x$ is the trivial forbidden thread at $x$, we know that there exists an arrow $\a$ leaving $x$ or an arrow $\b$ ending at $x$. Then put $\s(p_x)=\e(h_x) = -\s(\a)$, for the first case and put $\s(p_x)=\e(h_x) = -\e(\b)$, for the second case. Given a path $w$, denote by $\ell(w)$ its length, that is its number of arrows.

\begin{algo}[Avella - Alaminos and Geiss \cite {Avella-Alaminos2008}]\label{algoritmo de phi} \ Let $A=\K Q/I$ be a gentle bound quiver, for which all permitted and forbidden threads are determined.
\begin{enumerate}
 \item \begin{enumerate}
        \item Begin with a permitted thread $H_0$ of $A$,
        \item If $H_i$ is defined, let $P_i$ be the forbidden thread sharing target with $H_i$ and such that $\e(H_i)=-\e(P_i)$,
        \item Let $H_{i+1}$ be the permitted thread sharing source with $P_i$ and such that $\s(H_{i+1}) = -\s(P_i)$.
        \item[]The process stops if $H_n= H_0$ for some natural number $n$. In this case, let $m = \sum_{1\leqslant i \leqslant n} \ell(P_i)$
       \end{enumerate}
\item Repeat step 1 until all permitted threads of $A$ have been considered;
\item If there are (oriented) cycles $w$ with full relations, add a pair $(0,\ell(w))$ for each of those cycles;
\item Define \mor{\phi_A}{\mathbb{N} \times \mathbb{N} }{\mathbb{N}} by letting $\phi_A(n,m)$ be the number of times the pair $(n,m)$ appears in the algorithm.
\end{enumerate}
\end{algo}

Theorem A in \cite{Avella-Alaminos2008} asserts that $\phi$ is a derived invariant, but so far, it was only known to be a complete derived invariant in some particular cases,  see \cite[Theorem C]{Avella-Alaminos2008} or \cite{AA07}. Here, we show that it is not a complete derived invariant in general. Namely, we shall exhibit a class of gentle algebras for which $\phi$ is not a complete derived invariant.


\subsection{$m$-cluster tilted algebras} Let $H\simeq \K Q$ be an hereditary algebra. The derived category $\mathcal{D}^b(H)$ is triangulated, the translation functor, denoted by $[1]$, being induced from the shift of complexes. For an integer $n$, we denote by $[n]$ the composition of $[1]$ with itself $n$ times, thus $[1]^n = [n]$. In addition, $\mathcal{D}^b(H)$ has Auslander-Reiten triangles, and, as usual, the Auslander-Reiten translation is denoted by $\tau$.

Let $m$ be a natural number. The $m$-cluster category of $H$ is the quotient category $\mathcal{C}_m(H):=\mathcal{D}^b(H)/ \tau^{-1} [m]$ which carries a natural  triangulated structure, see \cite{K05}. Following \cite{Thomas2007} we consider \textit{$m$-cluster tilting objects}  in $\mathcal{C}_m(H)$ defined as objects satisfying the following conditions:
\begin{enumerate}
 \item $\Ho_{\mathcal{C}_m(H)}(T,X[i]) = 0$ for all $i\in\{1,2,\ldots,m\}$ if and only if $X\in \ts{add}\ T$,
 \item $\Ho_{\mathcal{C}_m(H)}(X,T[i]) = 0$ for all $i\in\{1,2,\ldots,m\}$ if and only if $X\in \ts{add}\ T$.
\end{enumerate}

The endomorphism algebras of such objects are called \emph{$m$-cluster tilted algebras of type $Q$}. In case $m=1$, this definition specializes to that of a cluster tilted algebra, a class intensively studied since its definition in \cite{BMR06}.

In \cite{ABCP09} it has been shown that cluster tilted algebras are gentle if and only if they are of type  $\mathbb{A}$ or $\tilde{\mathbb{A}}$. In \cite{Buan2008} Buan and Vatne gave the derived equivalence classification of cluster tilted algebras of type $\mathbb{A}$. They showed that two cluster tilted algebras of type $\mathbb{A}$ are derived equivalent if and only if their quivers have the same  number of 3-cycles with full relations and the same number of arrows.  Later, in \cite{Bastian09} the same work has been done for cluster tilted algebras of type $\tilde{\mathbb{A}}$. Moreover, in \cite{BB10}, the algebras that are derived equivalent to cluster tilted algebras of types $\mathbb{A}$ or $\tilde{\mathbb{A}}$ have been classified. In this classification, again, combinatorial data of the involved bound quiver is of central importance.

On the other hand, using arguments similar to those of \cite{ABCP09}, Murphy showed in \cite{Murphy2010} that $m$-cluster tilted algebras of type $\mathbb{A}$ are gentle and he described the connected components of $m$-cluster tilted algebras up to derived equivalence, a result analogous to that of \cite{Buan2008}. Later, a similar work has been done in \cite{Gubitosi} for  $m$-cluster tilted algebras of type $\mathbb{\widetilde{A}}$, where it is shown that $m$-cluster tilted algebras of type $\mathbb{\widetilde{A}}$ are gentle and their possible bound quivers are described. Moreover, in \cite{GubBust} the algebras that are derived equivalent to $m$-cluster tilted algebras of type  $\mathbb{A}$ have been classified. They are called branched algebras   \cite[Definition 4.3]{GubBust}.

\section{$m$-cluster tilted algebras of type $\mathbb{\widetilde{A}}$}


Given a bound quiver $(Q,I)$ and an integer $m$, a cycle is called \emph{ $m$-saturated} if it is an oriented cycle consisting of $m+2$ arrows such that the composition of any two consecutive arrows on this cycle belongs to $I$. Recall that two relations  $r$ and $r'$ in the bound quiver $(Q,I)$ are said to be   \textit{consecutive} if there is a walk  $v=wr=r'w'$ in  $(Q,I)$ such that $r$ and $r'$ point in the same direction and share an arrow.\\

\begin{defi}\cite[Definition 7.2]{Gubitosi}
 Let $\widetilde{\mathcal{C}}$ be a cycle without  relations (oriented or not) and fix an orientation of its arrows.  We say that an algebra  $A\cong \K Q/I$ is an  \textit{algebra with root $\widetilde{\mathcal{C}}$} if its bound quiver can be constructed  as follows:

\begin{enumerate}
\item We add to the  cycle $\widetilde{\mathcal{C}}$ gentle quivers in such a way that the final quiver remains gentle and connected. These added gentle quivers can only have  $m$-saturated cycles. We call these quivers  \textit{rays}.
\item We can add  relations to the  cycle $\widetilde{\mathcal{C}}$. If the cycle  $\widetilde{\mathcal{C}}$ is oriented then  we must add at least one relation.\\
\end{enumerate}

Also, we will refer to the  cycle $\widetilde{\mathcal{C}}$ as the  \textit{root cycle}.
\end{defi}

Let $\widetilde{\mathcal{C}}$ be a  cycle and $A$ an algebra with root  $\widetilde{\mathcal{C}}$. Each ray of $A$ can share with the cycle $\widetilde{\mathcal{C}}$  at most  $m+2$ vertices. If it shares just one vertex, this vertex is the \emph{union vertex} of the ray. If it shares more than one vertex, the ray and the cycle $\widetilde{\mathcal{C}}$ are connected through an  $m$-saturated cycle. For each union vertex there is at least one relation $\rho$ involving at least one arrow of $\widetilde{\mathcal{C}}$. If both arrows of $\rho$ belong to the root cycle, $\rho$ is called \textit{internal union relation} of the ray. If instead  just one arrow of $\rho$ belongs to the root cycle, $\rho$ is called \textit{external union relation} of the ray.

\begin{obs}\label{ciclo q comparte mas de una flecha entonces hay flechas - 1 relaciones en el sentido contrario}
\textnormal{Because of  \cite[Lemma 7.11]{Gubitosi} we know that if an  $m$-saturated cycle shares with the root cycle $k-1$ arrows counterclockwise oriented  and  $m+2-k+1$ arrows clockwise oriented,  then  one of the following holds:
\begin{itemize}
  \item [(a)] there is at least $k-2$ clockwise internal relations or at least  $m+2-k$ counterclockwise internal relations.
  \item [(b)] there is another  $m$-saturated cycle sharing with the root cycle $k-1$ arrows clockwise oriented  and  $m+2-k+1$ arrows counterclockwise oriented.
\end{itemize}}

\end{obs}

 For $m\geqslant 2$, an  $m$-cluster tilted algebras of type $\tilde{\mathbb{A}}$ (as is the case for type $\mathbb{A}$) does not have to be connected. However we will work with the connected component containing the non saturated cycle. Every other component corresponds  to an $m$-cluster tilted algebra of  type $\mathbb{A}$,  and consequently we know its derived class. See \cite{GubBust}.\\

\subsection*{Theorem }\cite[Theorem 7.16]{Gubitosi}\textit{ A connected algebra $A=\K Q/I$ is  a connected component of an $m$-cluster tilted algebra of type $\widetilde{\mathbb{A}}$ if and only if  $(Q,I)$ is a gentle bound quiver satisfying the following conditions: }

\textit{\begin{itemize}
\item [(a)] It  contains a non-saturated  cycle  $\widetilde{\mathcal{C}}$ in such a way that $A$ is an algebra with root $\widetilde{\mathcal{C}}$.
  \item [(b)] If it contains more cycles, then all of them are  $m$-saturated cycles.
  \item [(c)] Outside of an $m$-saturated cycle it can have at most  $m-1$ consecutive relations.
  \item [(d)] If  $\widetilde{\mathcal{C}}$ is an oriented  cycle, then it must have at least one internal relation.
  \item [(e)] If there are internal relations in the root cycle, then   the number  of clockwise oriented  relations   is equal modulo $m$ to the number of counterclockwise oriented.
  \end{itemize}}

As in the case of  $m$-cluster tilted algebras of type $\mathbb{A}$ with $(m\geqslant 2)$, the class of $m$-cluster tilted algebras of type $\tilde{\mathbb{A}}$ is not closed under derived equivalence; that is,  it is possible for an $m$-cluster tilted algebra of type $\tilde{\mathbb{A}}$  to be derived equivalent to an algebra which is not $m$-cluster tilted.

\begin{ejem}
Consider the following quiver $Q$:

$$\SelectTips{eu}{10}\xymatrix@C=.4pc@R=.2pc{ & & . \ar[ddrr]^{\b_1} & & &&&&&&&&&&& \\
&  & & & &&&&&&&&&&&&&&\\
. \ar[uurr]^{\b_4} & & & & .\ar[ddll]^{\b_2} & & . \ar[ll] & & .\ar[ll] & &  & &   & &   & &&& \\
&  & & & &&&&&&&&&&&&&&\\
&&  .\ar[uull]^{\b_3}  &  & &&&&&& . \ar[lluu] \ar[lldd] & & . \ar[ll]^{\a_3} && . \ar[ll]^{\a_2} && . \ar[ll]^{\a_1} \\
&  & & & &&&&&&&&&&&&&&\\
 & & & & .\ar[lluu] & & . \ar[ll] & & .\ar[ll]}$$

\

Let $I_1$ be the ideal generated by relations of the form  $\b_i\b_{i+1}$ for $i\in\{1,2,3,4\}$ where indices are to be read modulo $4$. Then $\K Q/I_1$ is a $2$-cluster tilted algebra of type $\tilde{\mathbb{A}}$. On the other hand let $I_2$ be $I_1$ plus the ideal generated by the  relations $\a_i \a_{i+1}$. The algebra $\K Q /I_2$ is not $2$-cluster tilted, but the two algebras are derived equivalent, as we shall see.
\end{ejem}

Our aim is to classify the algebras that are derived equivalent to $m$-cluster tilted algebras of type $\mathbb{\widetilde{A}}$. Since the derived class of $m$-cluster tilted algebras of type $\mathbb{A}$ is well understood \cite{GubBust}, it remains to find  the connected algebras that are derived equivalent to the  connected component having the root cycle.\\

\section{  $\tilde{\mathbb{A}}$-branched algebras}




Let $A$ be an algebra with root and let $\mathcal{S}$ be the set of all arrows in the quiver of $A$ not belonging to any   $m$-saturated cycle.

\begin{defi}
The   \textit{number of free clockwise arrows} in $\mathcal{S}$ is equal to the number of clockwise oriented arrows on the root cycle that are not  involved in any internal union relation   plus the number of clockwise  internal union relations  plus the number of arrows on the rays associated to clockwise union relations  (internal or external).\\

Dually, we define the number  of free counterclockwise arrows.
\end{defi}

We conjecture that the algebras that satisfy the following definition are the algebras that we are interested in.

\begin{defi}
We say that a connected algebra  $B=\K Q/I$  is \textit{$\tilde{\mathbb{A}}$}\textit{-branched}  if $B$ satisfies the following conditions:

\begin{itemize}
\item [(a)] There is a  cycle $\widetilde{\mathcal{C}}$ in $Q$ in such a way that $B$ is an algebra with root $\widetilde{\mathcal{C}}$.
\item [(b)] In the root cycle  the number  $r_h$  of clockwise oriented relations is the same modulo  $m$  that the number  $r_a$ of counterclockwise oriented relations.
\item [(c)] If $|r_h-r_a|=r=\alpha (m-1) + \beta$ (with $\beta< m-1$), then there must exist $r+1+\varepsilon$ free arrows  not belonging to any   $m$-saturated cycle on the clockwise sense if  $r_h>r_a$ or in the counterclockwise sense otherwise. Here,

    $$\varepsilon =\left\{
              \begin{array}{ll}
                \alpha - 1, & \hbox{\text{if }  $\beta=0$ ;} \\
                \alpha, & \hbox{\text{if}  $\beta\neq 0$ .}
              \end{array}
            \right.$$
\end{itemize}

\end{defi}


Let $\mathcal{R}$ be a ray of an $\tilde{\mathbb{A}}$-branched bound quiver. It follows from \cite[Proposition 6.5]{GubBust} that the algebra associated to $\mathcal{R}$ is tilting-cotilting equivalent to an algebra without relations. That is, an algebra whose quiver has $m$-saturated cycles possibly separated by arrows.

According to the next remarks we can move closer together the cycles and  assume  that any two of them  can be attached at the vertex that we choose.\\

\begin{obs}\label{forma normal de los rayos}
\textnormal{Applying the mutation  $\sigma'_{c}$ to the gentle bound quivers:}

$$\SelectTips{eu}{10}\xymatrix@C=.1pc@R=.3pc{     &  &  & . \ar@{.}@/^ .1pc/[rr] && b_{m+1}\ar[ddr] &  &&  & c_1 \ar@{.}@/^ .1pc/[rr] && . \ar[ddr] &  &&&&&&&&\\
&&&&&&&&&&&&  &&&&&&&&\\
z \ar@{.}[rr]  & & .\ar[ruu]    &&&  & b \ar[ddl]  &&  c \ar[uur] \ar[ll]  &&& & . \ar[ddl]   & & a \ar@{.}[ll] \\
&&&&&&&&&&&&&&&&&  &&& \\
  &  &  & . \ar@{.}@/_/[rr] \ar[uul] && b_1 &  &&  & c_{m+1} \ar@{.}@/_/[rr] \ar[uul] && .  &    &&&&&&&&&&}
  \SelectTips{eu}{10}\xymatrix@C=.1pc@R=.3pc{     &  &  & . \ar@{.}@/^ .1pc/[rr] && b_{m+1}\ar[ddr] &  &&  & c_1 \ar@{.}@/^ .1pc/[rr] && . \ar[ddr] &  &&&&&&&&\\
&&&&&&&&&&&&  &&&&&&&&\\
z \ar@{.}[rr]  & & .\ar[ruu]    &&&  & b \ar[ddl] \ar[rr] &&  c \ar[uur]   &&& & . \ar[ddl]   & & a \ar@{.}[ll] \\
&&&&&&&&&&&&&&&&&  &&& \\
  &  &  & . \ar@{.}@/_/[rr] \ar[uul] && b_1 &  &&  & c_{m+1} \ar@{.}@/_/[rr] \ar[uul] && .  &    &&&&&&&&&&}$$

 we obtain respectively the gentle bound quivers:

  $$\SelectTips{eu}{10}\xymatrix@C=.1pc@R=.3pc{     &  &  & . \ar@{.}@/^ .1pc/[rr] && b_{m+1}\ar[ddr] &  & c_{m+1} \ar@{.}@/^ .1pc/[rr] \ar[ddl] && .  &  &&&&&&&&\\
&&&&&&&&&&&&  &&&&&&\\
z \ar@{.}[rr]  & & .\ar[ruu]    &&&  & b \ar[ddl] \ar[ddr]  &&& & . \ar[uul]   & &  \\
&&&&&&&&&&&&&&&&&  & \\
  &  &  & . \ar@{.}@/_/[rr] \ar[uul] && b_1 &  &  c \ar@{.}@/_/[rr]  && . \ar[uur]  &    &&&&&&&&&& \\
  &&&&&&&&&&&&&&&&&\\
    & &    &&&  &      && c_1 \ar[uul] \ar@{.}[ddr]  &    & & \\
    &&&&&&&&&&&&&&&&&&\\
    & &    &&&  &      && &     a  & }
    \SelectTips{eu}{10}\xymatrix@C=.1pc@R=.3pc{     &  &  & . \ar@{.}@/^ .1pc/[rr] && b_{m+1}\ar[ddr] &  & c_{1} \ar@{.}@/^ .1pc/[rr]  && .\ar[ddr]   &  &&&&&&&&\\
&&&&&&&&&&&&  &&&&&&\\
z \ar@{.}[rr]  & & .\ar[ruu]    &&&  & b \ar[ddl] \ar[uur]  &&& & . \ar[ddl]   & &  \\
&&&&&&&&&&&&&&&&&  & \\
  &  &  & . \ar@{.}@/_/[rr] \ar[uul] && b_1 &  &  c \ar@{.}@/_/[rr] \ar[ddr]  \ar[uul]&& .  &    &&&&&&&&&& \\
  &&&&&&&&&&&&&&&&&\\
    & &    &&&  &      && c_{m+1} \ar@{.}[ddr]  &    & & \\
    &&&&&&&&&&&&&&&&&&\\
    & &    &&&  &      && &     a  & } $$

(where the vertex  $a$ is the union vertex of the ray). \\


\end{obs}


\begin{obs}\label{dos ciclos saturados se pueden pegar en el vertice que queramos}
\textnormal{Let $\mathcal{B}$ and $\mathcal{C}$ be two $m$-saturated cycles. Assume that they are attached at the vertex $c$. Applying the mutation $\sigma'_c$ we move the vertex $c$ one place in the direction of the  cycle $\mathcal{B}$.}
\end{obs}


\begin{center}
\begin{tabular}{lcccccl}
$\SelectTips{eu}{10}\xymatrix@C=.6pc@R=1pc{
   &&&\\
   b_{m+1} \ar@{.}@/^1.5pc/[rr] \ar[d]& \mathcal{B} & b_3 & \\
   b_{m+2} \ar[dr] & & b_2 \ar[u] &\\
   & c \ar[rd] \ar[ru] && \\
   c_{m+2} \ar[ru] &  & c_2 \ar[d] &\\
   c_{m+1} \ar[u]\ar@{.}@/_1.5pc/[rr] & \mathcal{C} & c_3\\
   &&&&
   }$
&&&
$\SelectTips{eu}{10}\xymatrix@C=.6pc@R=1pc{ & \\
&\\
&\\
\stackrel{\sigma'_c}{\longmapsto}\\
&\\
&\\
&}$
&&&
$\SelectTips{eu}{10}\xymatrix@C=.6pc@R=1pc{
   &&&\\
   b_{m+1} \ar@{.}@/^1.5pc/[rr] \ar[d]& \mathcal{B} & b_3 & c_3  \ar[r] & c_4 \ar@{.}@/^1.5pc/[dd] \\
   b_{m+2} \ar[dr] & & c \ar[u] \ar[ur] & & \mathcal{C} & \\
   & c_2 \ar[ru] && b_2 \ar[lu] & c_{m+2} \ar[l] \\
   &&&&
   }$

\end{tabular}
\end{center}

Therefore, we can assume that every ray has the following form:

$$\xymatrix@C=.3pc@R=.3pc{&&&   & .\ar@{.}@/^.1pc/[rrr] &&& .\ar[ddr] &  & . \ar@{.}@/^ .1pc/[rrr] &&& .\ar[ddr] &  &     &&&&  & . \ar@{.}@/^ .1pc/[rrr] &&& . \ar[ddr] &  &&&&&&&&&&\\
&&&  &&&&&&&&&&&&&&&&&&&&  &&&&&&&&&&\\
(*) &&&    z\ar[ruu] & &&& & . \ar[ddl] \ar[uur] &  &&&  & . \ar[ddl] \ar@{.}[rrrrr] & &&&&  .\ar[uur] & &&& & . \ar[ddl]  && . \ar@{-}[ll]  &&  . \ar@{-}[ll]\ar@{.}[rrrr] &&&& . \ar@{-}[ll]&& a \ar@{-}[ll]\\
&&&    &&&&&&&&&&&&&&&&&&&&  &&&&&&&&&& \\
&&&   & .\ar@{.}@/_/[rrr]  \ar[uul] &&& . &  & . \ar@{.}@/_/[rrr] \ar[uul] &&&  &  &  &&&&  & . \ar@{.}@/_/[rrr] \ar[uul] &&& .  &    &&&&&&&&&&}$$

\vspace*{.5cm}

where the vertex $a$ is the union vertex. Recall that here the orientation of the linear part is not fixed like in the normal form of  $m$-cluster tilted algebras of  type $\mathbb{A}$ given in \cite{Murphy2010}.

 In light of the preceding result, we introduce a somehow intermediate class of   $\tilde{\mathbb{A}}$-branched algebras.

 \begin{defi}
 An $\tilde{\mathbb{A}}$-branched algebra $B=\K Q/I$ is said to be \textit{solar}  if every ray of $(Q,I)$ has the form  $(*)$ above.

Given a ray  $\mathcal{R}$  in a solar algebra  let $a$ be its union vertex and  $\alpha$ the adjacent arrow. We say that the ray $\mathcal{R}$ is:
 \begin{itemize}
     \item [a)] \textit{incoming} to the root cycle if  $t(\alpha)=a$.
     \item [b)] \textit{outgoing} to the root cycle if  $s(\alpha)=a$.
   \end{itemize}

 \end{defi}

It follows from the previous remarks that every $\tilde{\mathbb{A}}$-branched algebra is   \BBeq to a solar algebra, then in the sequel we will work with the smaller  class of solar algebras.

%
%
%
%
%
%
%
%

%
%
%
%

%
%


\section{Normal forms}\label{normal forms}

It is our aim to show that every solar algebra with a given number $k$, of  $m$-saturated cycles is derived equivalent to one of  the following normal forms having the same number of $m$-saturated cycles.\\

\begin{defi}\label{non-oriented normal form}

For  $n_1,n_2\in \mathbb{N}^*$, $k_1,k_2\in \mathbb{N}$ and $r\in \mathbb{Z}$ we define the \emph{non-oriented normal form} $\widetilde{N}_{n_1,k_1,n_2,k_2,r}$  to be the following $\tilde{\mathbb{A}}$-branched algebra given as a quiver with relations:

$$\xymatrix@C=.4pc@R=.6pc{ &&& &&& &&& &&& &&& &&& &&& &&& & .\ar[dddl]_{\gamma_{1}^{m+1}} \ar@{.}@/^/[r]& . &\\
&&& &&& &&& &&& &&& &&& &&& &&& &&&\\
&&& &&& &&& &&& &&& &&& &&& &&& &&&\\
&&& .\ar[rrr]^{\a_1} & \ar@{.}@/^/[lld] & \ar@{.}@/_/[rr]&  . \ar[rrr]^{\a_2} && \ar@{.}@/_/[rr]& . \ar[rrr]^{\a_3} &&& \ar@{.}[rrr] &&\ar@{.}@/_/[rr]& . \ar[rrr]^{\a_r} &&& . \ar@{.}[rrr] &&& . \ar[rrr]^{\a_{r+s}} &&&  . \ar[rrr]_{\a_{r+s+1}} &&& .\ar[uuul]_{\gamma_{1}^1} \ar@{.}[rrrdd] &&&  && .\ar[ddll]_{\gamma_{k_2}^{m+1}} \ar@{.}@/^/[dr] &&& &&&\\
&&& &&& &&& &&& &&& &&& &&& &&& &&& &&& &&& . &&& &&& \\
 &&& &&& &&& &&& &&& &&& &&& &&& &&& &&&  . \ar[rrdd]_{\a_{r+s+k_2}}    \\
  &&& &&& &&& &&& &&& &&& &&& &&& &&& &&& &&& \\
 \scriptstyle{0} \ar[rrruuuu]^{\a_0} \ar[rrrdddd]_{\b_0} &&& &&& &&& &&& &&& &&& &&& &&& &&&  &&& && .\ar[uuur]_{\gamma_{k_2}^1} \ar[dddr]^{\d_{k_1}^1}&&& \\
 &&& &&& &&& &&& &&& &&& &&& &&& &&& &&& &&& \\
  &&& &&& &&& &&& &&& &&& &&& &&& &&& &&&   \ar[rruu]^{\b_{n+k_1}}  . &&& \\
  &&& &&& &&& &&& &&& &&& &&& &&& &&&  &&& &&& .  \\
 &&& .\ar[rrr]_{\b_1} &&& . \ar[rrr]_{\b_2} &&& . \ar[rrr]_{\b_3} &&& \ar@{.}[rrrrrrrrr] &&&   &&&   &&& . \ar[rrr]_{\b_n} &&&. \ar[rrr]^{\b_{n+1}} &&&  . \ar[dddl]^{\d_{1}^1} \ar@{.}[rrruu] &&& && . \ar[uull]^{\d_{k_1}^{m+1}} \ar@{.}@/_/[ur]  &&& &&&\\
&&& &&& &&& &&& &&& &&& &&& &&& &&&\\
&&& &&& &&& &&& &&& &&& &&& &&& &&&\\
 &&& &&& &&& &&& &&& &&& &&& &&& & . \ar[uuul]^{\d_{1}^{m+1}}\ar@{.}@/_/[r]& . &}$$

where every oriented cycle is $m$-saturated.

\begin{itemize}
  \item $n_1=n+1+k_1$ is the number of arrows counterclockwise oriented.
   \item $k_1$ is the number of counterclockwise  $m$-saturated  cycles.
  \item $n_2=r+s+1+k_2$  is the number of arrows clockwise oriented.
  \item $k_2$ is the number of clockwise  $m$-saturated  cycles.
  \item $r$ is the number of clockwise oriented relations if $r\geq 0$ or counterclockwise oriented if  $r<0$.
  \item $|r|\equiv 0$ modulo $m$.
  \item Si $r=\alpha (m-1) + \beta$ (with $\beta< m-1$), then  $s\geq r+1+\varepsilon$,  where

    $$\varepsilon =\left\{
              \begin{array}{ll}
                \alpha - 1, & \hbox{\text{if }  $\beta=0$ ;} \\
                \alpha, & \hbox{\text{if}  $\beta\neq 0$ .}
              \end{array}
            \right.$$
\end{itemize}

\end{defi}

\begin{defi}\label{forma normal orientada}

For $k,n,t \in \mathbb{N}^*$ and $n\geq 1$  we define the \emph{oriented normal form} $B_{k,n,t}$ to be the following $\tilde{\mathbb{A}}$-branched  algebra given as a quiver with relations:

$$\xymatrix@C=.4pc@R=.8pc{  &&& &&& &&& & .\ar[dddl]_{\gamma_{k}^{m+1}} \ar@{.}@/^/[r]& . &\\
 &&& &&& &&& &&&\\
 &&& &&& &&& &&&\\
 &&&  &&& \scriptstyle{n+t-1} \ar[rrr]^{\a_{n+t-1}} &&&  \scriptstyle{n+t} \ar[rrr]_{\a_{n+t}} &&& .\ar[uuul]_{\gamma_{k}^1} \ar@{.}@/^/[rrrrdd] &&&  && &  &&& &&&\\
 &&& &&& &&& &&& &&& &&&  & .\ar[dlll]_{\gamma_{2}^{m+1}} \ar@{.}@/^/[dr] &&& &&&  \\
  && &\scriptstyle{n+1}\ar@{.}@/^/[uurrr] &&& &&& &&& &&&  &  . \ar@/^.5pt/[rdd]_{}  &&& & . \\
   &&& &&& &&& &&& &&& &&& \\
  &&  &\scriptstyle{n} \ar[uu]^{\a_n} &&& &&& &&&  &&& && . \ar@/^.5pt/[ddl]_{\a_{n+t+k-1}} \ar[uurrr]_{\gamma_{2}^1} &&& \\
  && & \ar@{.}@/^/[ddr]&&& &&& &&& &&& &&& \\
  && &\scriptstyle{n-1}\ar[uu]^{\b_{n-1}} &&& &&  & & & && & && \scriptstyle{0} \ar[rrrd]_{\gamma_1^1}   \ar@/^.9pt/[ddll]^{\b_0}  &&& & . \ar[llluu]_{\gamma_1^{m+1}} \ar@{.}@/^/[dl]\\
 &&& &&&  &&  & &&&&&&& && &  . &&& &&&\\
   &&& & \scriptstyle{n-2} \ar[uul]^{\b_{n-2}}\ar@{.}@/_/[rrrrd]&&&& &&& &&\ar@{.}@/^/[rru]& \scriptstyle{1} \ar[llld]^{\b_1}&&&&&&\\
 && &&& &&& \scriptstyle{3} && \ar@{.}@/^/[rru]& \scriptstyle{2}  \ar[lll]^{\b_2}&&&&&&&&&&
}$$

\medskip

\begin{itemize}
   \item $k$ is the number of   $m$-saturated  cycles.
  \item $n-1$  is the number of consecutive relations  $\b_i\b_{i-1}$ ($i\in \{1, \cdots , n-1\}$).
  \item $n+t$ is the number of arrows not belonging to any  $m$-saturated cycle.
  \item Si $n-1=\alpha (m-1) + \beta$ (with $\beta< m-1$), then  $t\geq n+\varepsilon$,  where

    $$\varepsilon =\left\{
              \begin{array}{ll}
                \alpha - 1, & \hbox{\text{if }  $\beta=0$ ;} \\
                \alpha, & \hbox{\text{if}  $\beta\neq 0$ .}
              \end{array}
            \right.$$
\end{itemize}


\end{defi}

\vspace*{.5cm}

Given a pair $(a,b)\in \mathbb{N} \times \mathbb{N}$, denote by $(a,b)^\ast$ the characteristic function of the set  $\{(a,b)\} \subseteq \mathbb{N} \times \mathbb{N}$. The following will be useful.

\begin{prop}\label{phi de la forma normal A tilde}
Let $\widetilde{N}=\widetilde{N}_{n_1,k_1,n_2,k_2,r}$ be the non-oriented normal form. Then  the Avella-Alaminos and Geiss map, can be written as:

$$\phi_{\widetilde{N}}=((m-1)k_1+n_1+r,n_1-k_1)^*+ ((m-1)k_2+n_2-r,n_2-k_2)^*+ (k_1+k_2).(0,m+2)^*$$

\end{prop}

\begin{proof}

Let $(Q, I)$ be the bound quiver associated to the non-oriented normal form. Label the arrows of $Q$ as follows. Let $\b_0,\b_1,\ldots , \b_{n+k_1}$ be the counterclockwise arrows of the non-oriented cycle not belonging to any $m$-saturated cycle, where the source of  $\b_0$ is $0$, the unique source of the quiver,  and  $ \b_{n+k_1}\ldots \b_1 \b_0$ does not belong to $I$.  Furthermore, let  $\a_0,\a_1,\ldots , \a_{r+s+k_2}$ be the clockwise arrows, where the  source of $\a_0$ is also $0$,   $\a_{r+s+k_2}\ldots \a_{r+1} \a_r$ does not belong to  $I$ but $\a_{i}\a_{i-1}$ belongs to  $I$ for all $i\in\{1, \cdots, r\}$. For each $j\in \{1,2\ldots ,k_2\}$ let $\gamma_j^1, \gamma_j^2,\ldots , \gamma_j^{m+1}, \a_j^{r+s+j}$ be the  $m+2$  arrows of the  $j^{\text{-th}}$ cycle clockwise oriented and for each   $j\in \{1,2\ldots ,k_1\}$, let $\d_j^1, \d_j^2,\ldots , \d_j^{m+1}, \b_j^{n+j}$ be the  $m+2$ arrows of the   $j^{\text{-th}}$ cycle counterclockwise oriented. See definition \ref{non-oriented normal form}.

The term $(k_1+k_2) \cdot (0,m+2)^\ast$ comes from step $(3)$ of the algorithm \ref{algoritmo de phi}.

We start the algorithm with the permitted thread  $H_0=h_{t(\b_0)}$, then $\Pi_0=\b_0$ and the permitted thread   sharing source with the forbidden thread $\Pi_0$ is $H_1=\a_0$.
The algorithm can be summarized in the following table:

$$\begin{array}{lcl}
H_0 = h_{t(\b_0)}  &   &  \Pi_0 = \b_0 \\
H_1=\a_0          &   &  \Pi_1 = p_{t(\a_0)}\\
H_2=\a_1          &   &  \Pi_2 = p_{t(\a_1)}\\
   \vdots &             &            \vdots  \\
H_r = \a_{r-1}     &  &  \Pi_r = p_{t(\a_{r-1})}\\
H_{r+1} = \d_{k_1}^1\cdot\a_{r+s+k_2}\cdots \a_{r}  &  & \Pi_{r+1}= p_{t(\d_{k_1}^1)}\\
H_{r+2} = \d_{k_1^2} &  & \Pi_{r+2}= p_{t(\d_{k_1}^2)}  \\
   \vdots &             &             \vdots \\
H_{r+m}= \d_{k_1}^m &  & \Pi_{r+m} =  p_{t(\d_{k_1}^m)}\\
H_ {r+m+1}=  \d_{k_1-1}^1\cdot\d_1^{m+1} &  &  \Pi_{r+m+1}=p_{t(\d_{k_1-1}^1)} \\
H_ {r+m+2}=\d_{k_1-1}^2  &  &  \Pi_{r+m+2}= p_{t(\d_{k_1-1}^2)}\\
\vdots &             &            \vdots  \\
H_ {r+2m} =\d_{k_1-1}^m    &  & \Pi_{r+2m}=p_{t(d_{k_1-1}^m)} \\
\vdots &             &            \vdots  \\
H_ {r+(k_1-1)m+1}=\d_{1}^1 \cdot \d_2^{m+1} && \Pi_{r+(k_1-1)m+1} = p_{t(\d_{1}^1)}\\
H_ {r+(k_1-1)m+2}=\d_1^2  && \Pi_ {r+(k_1-1)m+2}=p_{t(\d_1^2)}\\
\vdots &             &            \vdots  \\
H_{r+k_1m} = \d_1^m && \Pi_{r+k_1m} =p_{t(\d_1^m)}\\
H_{r+k_1m+1}= \d_1^{m+1} && \Pi_{r+k_1m+1} =\b_n\\
H_ {r+k_1m+2}=h_{s(\b_n)} && \Pi_{r+k_1m+2} =\b_{n-1}\\
\vdots &             &            \vdots  \\
H_ {r+k_1m+n}= h_{s(\b_2)} && \Pi_{r+k_1m+n} = \b_1\\
H_ {r+k_1m+n+1} = h_{s(\b_1)} = H_0 && \\
\end{array}$$

The only forbidden threads of non-zero length are the arrows $\b_i$ with  $0\in \{1,\cdots, n\}$. Then, we get the pair $(r+k_1m+n+1,n+1)$, and since $n+1=n_1-k_1$, we obtain the pair $(r+(m-1)k_1+n_1,n_1-k_1)$.

Since not all permitted threads have been considered, we continue the algorithm with the permitted thread  $H_0=\gamma_{k_2}^1\b_{n+k_1}\cdots \b_0$. The following table summarizes this part of the algorithm.

$$\begin{array}{lcl}
H_0 = \gamma_{k_2}^1\b_{n+k_1}\cdots \b_0  &   &  \Pi_0 = p_{t(\gamma_{k_2}^1)}\\
H_1=\gamma_{k_2}^2          &   &  \Pi_1 = p_{t(\gamma_{k_2}^2)}\\
   \vdots &             &            \vdots  \\
H_{m-1} = \gamma_{k_2}^m    &  &  \Pi_{m-1} = p_{t(\gamma_{k_2}^m)}\\
H_{m} = \gamma_{k_2-1}^1\cdot\gamma_{k_2}^{m+1} &  & \Pi_{m}= p_{t(\gamma_{k_2-1}^1)}\\
H_{m+1} = \gamma_{k_2-1}^2 &  & \Pi_{m+1}= p_{t(\gamma_{k_2-1}^2)}  \\
   \vdots &             &             \vdots \\
H_{m+(m-1)}= \gamma_{k_2-1}^m &  & \Pi_{m+m-1} =  p_{t(\d_{k_1}^m)}\\
H_{2m}=  \gamma_{k_1-1}^1\cdot\gamma_1^{m+1} &  &  \Pi_{2m}=p_{t(\gamma_{k_1-1}^1)} \\
\vdots &             &            \vdots  \\
H_ {(k_2-1)m}=\gamma_{1}^1 \cdot \gamma_2^{m+1} && \Pi_{(k_2-1)m} = p_{t(\gamma_{1}^1)}\\
H_ {(k_2-1)m+1}=\gamma_1^2  && \Pi_ {(k_2-1)m+1}=p_{t(\gamma_1^2)}\\
\vdots &             &            \vdots  \\
H_{(k_2-1)m+(m-1)} = \gamma_1^m && \Pi_{k_2m-1} = \a_{r+s}\\
H_ {k_2m}=h_{s(\a_{r+s})} && \Pi_{k_2m} =\a_{r+s-1}\\
\vdots &             &            \vdots  \\
H_ {k_2m+s-1}= h_{s(\a_{r+2})} && \Pi_{k_2m+s-1} = \a_{r+1}\\
H_ {k_2m+s}= h_{s(\a_{r+1})} && \Pi_{k_2m+s} = \a_r\cdots \a_0\\
H_ {k_2m+s+1} =  H_0 && \\
\end{array}$$

\medskip

The only forbidden threads of non-zero length are the $s$ arrows $\a_{r+i}$ with  $i\in \{1,\cdots, s\}$  of length $1$ and the path $\a_r\cdots\a_0$ of length $r+1$. Then we get the pair $(mk_2+s+1,s+r+1)$. Using the equation $s+r+1=n_2-k_2$ we rewrite this pair as $((m-1)k_2+n_2-r,n_2-k_2)$.

Since all permitted threads have been considered, the algorithm  ends.
\end{proof}

\begin{coro}\label{si dos formas normales son der eq entonces r y r' son cong}
Let $N=\widetilde{N}_{n_1,k_1,n_2,k_2,r}$ and $N'=\widetilde{N}_{n_1',k_1',n_2',k_2',r'}$ be two derived equivalent normal forms. Then $|r|\equiv |r'|$ (mod $m$).
\end{coro}

\begin{proof}
Since $N$ and  $N'$ are derived equivalent we have $\phi_N=\phi_{N'}$.
 From $\phi_N$ we get, in particular,  the pairs $((m-1)k_1+n_1+r,n_1-k_1)=(mk_1+r+n+1,n+1)$  and $((m-1)k_2+n_2-r,n_2-k_2)=(mk_2+s+1,r+s+1)$ and, in consequence, the numbers $mk_1+r+n+1-(n+1)=mk_1+r$ and $mk_2+s+1-(r+s+1)=mk_2-r$. Doing the same with $\phi_{N'}$ we get the numbers $mk_1'+r'$ and  $mk_2'-r'$. Considering these numbers modulo $m$ we get the set $\{r,-r\}$ for $N$ and the set $\{r',-r'\}$ for $N'$.  We conclude that $|r|\equiv |r'|$ (mod $m$).
\end{proof}

\begin{prop}\label{phi de la forma normal orientada}
Let $B_{k,n,t}$  be the oriented normal form. Then,

 $$\phi_{B_{k,n,t}} = (n-1,0)^* + (t+mk+1,n+t)^* + k\cdot (0,m+2)^*$$

\end{prop}

\begin{proof}
Let $(Q, I)$ be the bound quiver associated to the oriented normal form. Label the arrows of $Q$ as follows.  For each $j\in \{1,2\ldots ,k\}$ let $\gamma_j^1, \gamma_j^2,\ldots , \gamma_j^{m+1}, \a_{n+t+k-j}$ be the  $m+2$  arrows of the  $j^{\text{-th}}$ $m$-saturated cycle.  Let  $\b_0,\b_1,\ldots , \b_{n-1}$ and  $\a_{n},\ldots ,\a_{n+t-1} $ be the  arrows not belonging to any $m$-saturated cycle, where  $\b_{i}\b_{i-1}$ belongs to  $I$ for all $i\in\{1, \cdots, n-1\}$ and  $\a_{n}\ldots \a_{n+t-1} $   does not belong to  $I$. See definition \ref{forma normal orientada}.

We begin the algorithm with the permitted threads  $H_0=\b_0 \a_{n+t+k-1}\cdots \a_n\b_{n-1}$. Then,

$$\begin{array}{lcl}
H_0 = \b_0 \a_{n+t+k-1}\cdots \a_n\b_{n-1}  &   &  \Pi_0 = p_{1}\\
H_1=\b_1          &   &  \Pi_1 = p_2\\
   \vdots &             &            \vdots  \\
H_{n-2} = \b_{n-2}    &  &  \Pi_{n-2} = p_{n-1}\\
H_{n-1} =  H_0 && \\
\end{array}$$

and we get the pair $(n-1,0)$.

We continue the algorithm with the permitted thread $H_0=h_n$. Then,

$$\begin{array}{lcl}
H_0 = h_{n}  &   &  \Pi_0 = \b_{n-1} \cdots \b_0 \\
H_1=\gamma_1^1          &   &  \Pi_1 = p_{t(\gamma_1^1)}\\
H_2=\gamma_1^2         &   &  \Pi_2 = p_{t(\gamma_1^2)}\\
   \vdots &             &            \vdots  \\
H_m = \gamma_1^m        &  &  \Pi_r = p_{t(\gamma_1^m )}\\
H_{m+1} = \gamma_{2}^1\cdot\gamma_1^{m+1}\cdots \a_{r}  &  & \Pi_{m+1}= p_{t(\gamma_{2}^1)}\\
H_{m+2} = \gamma_2^2    &  & \Pi_{m+2}= p_{t(\gamma_{2}^2)}  \\
   \vdots &             &             \vdots \\
H_{2m}= \gamma_2^m      &  & \Pi_{2m} =  p_{t(\gamma_{2}^m)}\\
H_ {2m+1}=  \gamma_{3}^1\cdot\gamma_2^{m+1} &  &  \Pi_{2m+1}=p_{t(\gamma_{3}^1)} \\
\vdots &             &            \vdots  \\
H_ {(k-1)m+1} =\gamma_{k}^1\cdot \gamma_{k-1}^{m+1}    &  & \Pi_{(k-1)m+1}=p_{t(\gamma_{k}^1)} \\
H_ {(k-1)m+2}=\gamma_k^2  && \Pi_ {(k-1)m+2}=p_{t(\gamma_k^2 )}\\
\vdots &             &            \vdots  \\
H_{km} = \gamma_{k}^{m} && \Pi_{km} =p_{t(\gamma_{k}^{m})}\\
H_{km+1}= \gamma_k^{m+1} && \Pi_{km+1} =\a_{n+t-1}\\
H_ {km+2}=h_{n+t-1} && \Pi_{km+2} =\a_{n+t-2}\\
\vdots &             &            \vdots  \\
H_ {km+t}= h_{n+1} && \Pi_{km+t} = \a_n\\
H_ {km+t+1} = h_{n} = H_0 && \\
\end{array}$$

The only forbidden threads of non-zero  length are the   $t$ arrows $\a_{n+i}$ with $i\in \{0,\cdots, t-1\}$  of length $1$ and the path $\b_{n-1}\cdots\b_0$ of length  $n$. Then we obtain the pair $(t+mk+1, n+t)$.

Since all permitted threads have been considered, this part of the algoritm is over.
The term  $k \cdot (0,m+2)^\ast$ comes from step $(3)$ of the algorithm.
\end{proof}

\begin{coro}
The algebras associated to the normal forms  $B_{k,n,t}$ and   $\widetilde{N}_{n_1,k_1,n_2,k_2,r}$ are not derived equivalent.
\end{coro}

\begin{proof}
 If $B_{k,n,t}$ is derived equivalent to  $N=\widetilde{N}_{n_1,k_1,n_2,k_2,r}$ the functions  $\phi_{N}$ and $\phi_{B_{k,n,t}}$ must be identical. We will see that this is impossible. Assume that

$$\begin{array}{l}
 ((m-1)k_1+n_1+r, n_1-k_1)=(t+mk+1,n+t)  \\
  ((m-1)k_2+n_2-r, n_2-k_2)=(n-1,0)\\
  k\cdot (0,m+2)=(k_1+k_2)\cdot (0,m+2)
\end{array}$$

Moreover, since the number of arrows is also a derived invariant, we have the equation
$$n_1+n_2=n+t+k(m+2)$$

Then, we have  $(m-1)k_1+n_1+r=t+mk+1=n_1+n_2-n-2k+1$ and $r=n_2-n+1-2k_2-(m+1)k_1$ which implies that  $n-1=(m-1)k_2+n_2-r=(m+1)(k_1+k_2)+n-1$ and  $(m+1)(k_1+k_2)$ has to be zero. Then, $k_1=k_2=k=0$ and since $n_2-k_2=0$ we obtain  $n_2=0$, which is absurd.

If instead we assume that

$$\begin{array}{l}
((m-1)k_1+n_1+r, n_1-k_1)=(n-1,0)\\
((m-1)k_2+n_2-r, n_2-k_2)=(t+mk+1,n+t) \\
k\cdot (0,m+2)=(k_1+k_2)\cdot (0,m+2)
\end{array}$$

an analogous computation gives  the absurdity  $n_1=0$.
\end{proof}

\section{Toward the reduction to Normal Form}

The idea is to apply a sequence of mutations that remove the relations of a solar algebra that lie outside the $m$-saturated cycles. In fact, every external union relation can be removed but the internal union relations can only  be removed by pairs, one in the clockwise sense with other in the counterclockwise sense, as in the case of iterated tilted algebras of type $\tilde{\mathbb{A}}$, see \cite{AS87}. At the same time  we want to move every ray of a solar algebra into the root cycle.

In the sequel we adopt the following convention concerning decorations on the names of vertices:
$\xymatrix@C=.6pc@R=.25pc{& \\
a \ar@{-}@/^/[ru]  & }$ , $\xymatrix@C=.6pc@R=.25pc{
a \ar@{-}@/_/[rd]  &  \\&}$ , $\SelectTips{eu}{10}\xymatrix@C=.5pc@R=.4pc{ \\
&  a  \ar@{-}@/_/[ul]}$  or $\SelectTips{eu}{10}\xymatrix@C=.5pc@R=.4pc{
&  a  \ar@{-}@/^/[dl]\\ &}$
  means that the vertex $a$ belongs to the root cycle.

Let $\rho$ be an (external or internal) union  relation and  $\mathcal{R}$ the corresponding ray.  The following  lemma  (and its dual) allow  to assume  that the linear part of the ray is all oriented in the sense of the arrow adjacent to the root cycle.

\begin{lema}\label{la flecha en el otro sentido}\label{rel interna meter flecha en sentido contrario}
Let $\gamma$ be the first arrow of the ray oriented  in the opposite direction of the root cycle  adjacent arrow.
\begin{itemize}
  \item [(a)] Assume that $\mathcal{R}$ is outgoing to the root cycle and $\rho$ is an external  union  relation. Then, the gentle bound quiver

 $$\xymatrix@C=.6pc@R=.4pc{&&&&&&&&&&&& &\\
 &&&&&&&&&&&& d \ar[dd] \ar@{-}@/^/[ru]  & \\
 &&&&&&&&&&&& &\\
 \cdots & &. \ar[ll] \ar[rr]^{\gamma} && c_1 && c_2 \ar[ll]  && c_n \ar@{.}[ll] && c  \ar[ll] & \ar@{.}@/^/[ru]^{\rho}& b\ar[ll] \ar[dd] &   \\
 &&&&&&&&&&&& & \\
 &&&&&&&&&&&& a \ar@{-}@/_/[rd]& \\
  &&&&&&&&&&&& &}$$

is tilting-cotilting equivalent to the gentle bound quiver

 $$\xymatrix@C=.6pc@R=.4pc{&&&&&&&&&&&& &\\
 &&&&&&&&&&&& d \ar[dd] \ar@{-}@/^/[ru]  & \\
 &&&&&&&&&&&& &\\
 \cdots & &. \ar[ll]  && c_1 \ar[ll]  && c_2 \ar[ll] \ar@{.}[rr] && c_{n-1} \ar@{.}[ll] && c_n  \ar[ll] & \ar@{.}@/^/[ru]^{\rho}& c\ar[ll] \ar[dd]^{\gamma'} &   \\
 &&&&&&&&&&&& & \\
 &&&&&&&&&&&& b\ar[dd]\\
 &&&&&&&&&&&&&\\
 &&&&&&&&&&&& a \ar@{-}@/_/[rd]& \\
  &&&&&&&&&&&& &}$$

  \item [(b)]Assume that $\mathcal{R}$ is incoming to the  root cycle and $\rho$ is an internal  union  relation. Then, the  gentle bound quiver

$$\SelectTips{eu}{10}\xymatrix@C=.5pc@R=.4pc{ \\
&  x \ar[dd] \ar@{-}@/_/[ul]\\
& \ar@{.}@/_/[dd] \\
& c\ar[dd] & & c_1  \ar[ll] && c_2 \ar[ll] &&  c_{n-1} \ar@{.}[ll] && c_n \ar[ll] \ar[rr]^{\gamma} && . && \cdots \ar[ll] \\
& \\
& z \ar@{-}@/^/[dl]\\
&  }$$

\medskip

is tilting-cotilting equivalent to the gentle bound quiver

$$\SelectTips{eu}{10}\xymatrix@C=.5pc@R=.4pc{ \\
&  x \ar[dd] \ar@{-}@/_/[ul]\\
& \ar@{.}@/_9pt/[dd] \\
& c_1\ar[dd] & & c_2  \ar[ll] && c_{n-1} \ar@{.}[ll]  && c_n \ar[ll]  && . \ar[ll] && \cdots \ar[ll] \\
& \\
& c\\
&\\
& z \ar[uu] \ar@{-}@/^/[dl]\\
&  }$$
\end{itemize}

\end{lema}

\begin{proof}\

\begin{itemize}
  \item [(a)] Apply the sequence of mutations $\sigma_{c_n}\cdots \sigma_{c_1}$ to move   $n$ places to the right the arrow $\gamma$. Finally  apply the mutation $\sigma_c$ to take the arrow $\gamma $ onto the root cycle.
  \item [(b)] Apply the sequence $\sigma'_{c_1}\cdots \sigma'_{c_{n}}$ to move $n$ places to the left the arrow  $\gamma$. Finally applying $\sigma'_c$ we take $\gamma $ inside the root cycle.
\end{itemize}
\end{proof}

%


\subsection{External union relations}

We are now able to describe the sequence of mutations that allows to remove the external union relations.  We start with an external union relation whose corresponding ray is outgoing to the root cycle. If the ray is incoming then the process is dual.

\begin{lema}\label{relacion externa sin ciclos saturados}
The gentle bound quiver

 $$\xymatrix@C=.6pc@R=.4pc{&&&&&&&&&& &\\
 &&&&&&&&&& a \ar[dd] \ar@{-}@/^/[ru]  & \\
 &&&&&&&&&& &\\
 && c_1 && c_2 \ar[ll]  && c_{n-1} \ar@{.}[ll] && c_n  \ar[ll] & \ar@{.}@/^/[ru]^{\rho}& b\ar[ll] \ar[dd] &   \\
 &&&&&&&&&& & \\
 &&&&&&&&&& b_1 \ar@{-}@/_/[rd]& \\
  &&&&&&&&&& &}$$

is tilting-cotilting equivalent to the gentle bound quiver

$$\xymatrix@C=.6pc@R=.25pc{& \\
a \ar[dd] \ar@{-}@/^/[ru]  & \\
&\\
c_n \ar[dd]  & \\
 &\\
 \vdots \ar[dd] &\\
 &\\
c_2 \ar[dd] &\\
&\\
c_1 \ar[dd] & \\
&\\
b \ar[dd]& \\
& \\
b_1 \ar@{-}@/_/[rd] & \\
  &}$$

\end{lema}

\begin{proof}
Apply the the sequence  $\sigma_{c_2}\cdots \sigma_{c_n}$ to move  the arrows not involved in $\rho$ into the root cycle. Finally, to remove the relation  $\rho$ and take the involved arrow inside the root cycle   apply $\sigma_{c_1}$.
\end{proof}

If we are given a solar algebra $A=\K Q/I$ with an external union relation $\rho$  whose corresponding ray $\mathcal{R}$ does not have $m$-saturated cycles, the previous lemmas show how to obtain an algebra $A'=\K Q'/I'$ which is tilting-cotilting equivalent to $A$. Moreover in $(Q',I')$ there is one relation less than in $(Q,I)$, namely the relation $\rho$ and the ray $\mathcal{R}$ was moved into the root cycle. Thus, it remains to see how to remove an external union relation such that the corresponding ray has $m$-saturated cycles.


In the sequel every  oriented cycle is an $m$-saturated cycle.

\begin{lema}\label{si tiene ciclos el rayo}

The  gentle bound quiver

$$\xymatrix@C=.6pc@R=.4pc{&&&&& &\\
 &&&&& a \ar[dd] \ar@{-}@/^/[ru]  & \\
 &c_m \ar[r] &c_{m+1} \ar[dr]&&& &\\
\ar@{.}@/^/[ur]  &&& c  \ar[dl] & \ar@{.}@/^/[ru]^{\rho}& b\ar[ll] \ar@{-}[dd] &   \\
  &c_2 \ar@{.}@/^/[lu] & c_1 \ar[l] & && & \\
 &&&&& b_1 \ar@{-}@/_/[rd]& \\
  &&&&& &}$$

 is tilting-cotilting equivalent to the gentle bound quiver

 $$\xymatrix@C=.6pc@R=.4pc{&& &&\\
 &&& a \ar[dd] \ar@{-}@/^/[ru]  & \\
 &&&&\\
 &&& c_{m+1} \ar[dd]  & \\
&c_{m-1} \ar[r] & c_m \ar[dr] && \\
 &&&c \ar[dd] &\\
 \ar@{.}@/^/[uru]&&&& \\
   & & & b  \ar[dl] \ar@{-}[dd] &   \\
  &c_2\ar@{.}@/^/[ulu] & c_1 \ar[l]&&  \\
 && & b_1 \ar@{-}@/_/[rd]& \\
  &&&& }$$

\end{lema}

\begin{proof}
Apply the sequence $\sigma_{c_{m+1}} \sigma_{c}$.
\end{proof}

Thus, until now, we know how to eliminate the external union relations and move the corresponding rays into the root cycle using mutations.


\subsection{Internal union relations}

%

Let $\rho$ be an  internal union  relation and  $\mathcal{R}$ the corresponding ray.
Now, we want to see how to move the ray $\mathcal{R}$ into the root cycle. We start with $\mathcal{R}$ having a non-empty linear part. In light of  lemma \ref{rel interna meter flecha en sentido contrario}, we can assume that all arrows in this linear part are oriented in  the same direction.

\begin{lema}\label{rel interna meter rayo}

The  gentle bound quiver

$$\SelectTips{eu}{10}\xymatrix@C=.2pc@R=.2pc{ &&&&&&&&&&&&&&&&\\
& x \ar@{-}@/_/[ul] \ar[ddd]&&&&&&&&&&&&&&\\
&&&& &  &&&&&&&&          & c_1 \ar@{.}@/^.1pc/[rrr] &&& .\ar[ddr] &  & . \ar@{.}@/^ .1pc/[rrr] &&& .\ar[ddr] &  &  &&&&  &\ar@{.}@/^ .1pc/[rrr] &&& . \ar[ddr] & \\
   & \ar@{.}@/_9pt/[dd] &&&&&&&&              &&&&&&&&&&&&&&&&&&&&  \\
     & a_1 \ar[ddd] &&& a_2 \ar[lll] &&& . \ar[lll] &&& . \ar@{.}[lll] &&& a_n \ar[lll]        \ar[ruu] & & \mathcal{C}_1 && & . \ar[ddl] \ar[uur] &  &  \mathcal{C}_2&&  & . \ar[ddl] \ar@{.}[rrrrr] & &&&&  .\ar[uur] & &  \mathcal{C}_r && & . \ar[ddl] \\
   &&&& &&&&&&&&&              &&&&&&&&&&&&&&&&&&&&   \\
    &&&& &  &&&&&&&&       & c_{m+1}\ar@{.}@/_/[rrr]  \ar[uul] &&& . &  & . \ar@{.}@/_/[rrr] \ar[uul] &&&  &  &  &&&&  & . \ar@{.}@/_/[rrr] \ar[uul] &&& .  & \\
   & z \ar@{-}@/^/[dl] &&&&&&&&&&&&&&&  \\
   &&&&&&&&&&&&&&&&&&&&&&&&& }$$

\vspace*{.5cm}


is tilting-cotilting equivalent to the gentle bound quiver

$$\SelectTips{eu}{10}\xymatrix@C=.3pc@R=.2pc{ &&&&&&&&&&&&&&&&&&&\\
                           & x \ar@{-}@/_/[ul] &&&&&&&&&&&\\
                           &&&&&&&&&&&&&&&&&\\
                           & a_1 \ar[dd] \ar[uu] &&&&&&&&&&\\
                           &&&&&&&&&&&&&&&&&\\
                           & a_2 \ar@{.}[dd] &&&&&&&&&&&&&&&&&\\
                           &&&&&&&&&&&&&&&&&\\
                           & .\ar[dd] & &&& &&&&&&&&&&&&&&&\\
                           &  &  c_1 \ar@{.}@/^/[rrr]&  && .\ar[ddr]& & . \ar@{.}@/^/[rrr]  &&& . \ar[ddr] &&  &&& . \ar@{.}@/^/[rrr] &&& . \ar[ddr]&&&&& \\
                           &  a_{n-1} \ar[ru] & & &  && &&&&& \\
                           &  && \mathcal{C}_1  && & .  \ar[uur] \ar[ddl] &&& \mathcal{C}_2 && . \ar@{.}[rrr] \ar[ddl] & && . \ar[ruu]  & &&\mathcal{C}_r && . \ar[ddl]   \\
                           &  a_n \ar[dd]   \ar[uu]  &&&&&&&&&&\\
                           & \ar@{.}@/_10pt/[dd]  & c_m \ar@{.}@/_/[rrr]\ar[ul]&&& . && . \ar[uul]\ar@{.}@/_/[rrr] &&& . & &&&  & .\ar@{.}@/_/[rrr] \ar[uul] &&& . &&\\
                           & c_{m+1} \ar[dd] &&&&&&&&&&&&&&&&&&\\
                           &&&&&&&&&&&&&&&&&&& \\
                           & z \ar@{-}@/^/[dl] &&&&&&&&&&&&&&&\\
                           &&&&&&&&&&&&&&&&&&&}$$

\end{lema}

\begin{proof}

Use the sequence $\sigma_{a_{n-1}}\cdots \sigma_{a_1}$ to move inside the root cycle the linear part of the ray. Then  apply $\sigma_{a_n}$ to move the $m$-saturated cycle $\mathcal{C}_1 $ inside the root cycle.
\end{proof}

\begin{obs}\label{rel interna meter ciclos}
Observe that in light of  lemma \ref{dos ciclos saturados se pueden pegar en el vertice que queramos} we can assume that the cycles $\mathcal{C}_1$ and $\mathcal{C}_{2}$ are attached at the vertex $c_m$. Then we can  move $\mathcal{C}_{2}$ in such a way that it shares one arrow with the root cycle by applying the mutation $\sigma'_{c_m}$:

\begin{minipage}[c]{0.97\linewidth}
\begin{minipage}{0.43\linewidth}
$$\SelectTips{eu}{10}\xymatrix@C=.2pc@R=.1pc{ &&&&&&&&&&&&&&&&&&&\\
                           & x \ar@{-}@/_/[ul] &&&&&&&&&&&\\
                           &&&&&&&&&&&&&&&&&\\
                           &&&&&&&&&&&&&&&&&\\
                           & a_1 \ar@{.}[ddd] \ar[uuu] &&&&&&&&&&\\
                           &&&&&&&&&&&&&&&&&\\
                           &&&&&&&&&&&&&&&&&\\
                           & .\ar[ddd] & &&& &&&&&&&&&&&&&&&\\
                           &&&  & &&  &&&&&&&&&&&&&&&\\
                           &  &  &&& & & &&&&& \\
                           &  a_{n-1} \ar[rr] & & c_1 \ar@{.}@/^1pc/[rrrdd] &&&&&  \\
                           &  &&   && &  &&&    \\
                           &&&&&&c_{m-1}  \ar[dlll] &&&&&&&&&&& &&& \\
                           & a_n   \ar[uuu] \ar[dd] && c_m \ar[ll]\ar[ddd]  &&&  &&  &&\\
                           &   \ar@{.}@/_9pt/[dd]   & & &&& b_{m+2} \ar[ulll] &&&&&&& \\
                           & c_{m+1} \ar[dd] &&&&\mathcal{C}_2 &&&&&&&&&&&&\\
                           &&& b_2  \ar[rrr] &&& b_3\ar@{.}@/_1.5pc/[uu]&&&&&&&&&&&&&& \\
                           & z \ar@{-}@/^/[dl] &&&&&&&&&&&&&&&\\
                           &&&&&&&&&&&&&&&&&&&}$$
\end{minipage}
\hfill
\begin{minipage}{0.1\linewidth}
$\substack{\stackrel{\sigma'_{c_m}}{\longmapsto}}$
\end{minipage}
\hfill
\begin{minipage}{0.43\linewidth}
$$\SelectTips{eu}{10}\xymatrix@C=.1pc@R=.05pc{ &&&&&&&&&&&&&&&&&&&\\
                           & x \ar@{-}@/_/[ul] &&&&&&&&&&&\\
                           &&&&&&&&&&&&&&&&&\\
                           &&&&&&&&&&&&&&&&&\\
                           & a_1 \ar@{.}[ddd] \ar[uuu] &&&&&&&&&&\\
                           &&&&&&&&&&&&&&&&&\\
                           &&&&&&&&&&&&&&&&&\\
                           & .\ar[ddd] & &&& &&&&&&&&&&&&&&&\\
                           &&& & &&  &&&&&&&&&&&&&&&\\
                           &  &  &&& & & &&&&& \\
                           &  a_{n-1} \ar[rr] & & c_1 \ar@{.}@/^1pc/[rrrdd] &&&&&  \\
                           &  &&   && &  &&&    \\
                           &&&&&&c_{m-1}  \ar[dlll] &&&&&&&&&&& &&& \\
                           & c_m   \ar[uuu] \ar[drr]   &&b_2 \ar[ll]  &&&  &&  &&\\
                           && & b_3  &&&  &&& &&&&&&&&&              &&&&  \\
                            &     & \mathcal{C}_2 &&&  &&&&&&& \\
                           & a_n \ar[ddd] \ar[uuu] && b_{m+2}\ar@{.}@/_1.5pc/[uu] \ar[ll]  &&&&&&&&\\
                           &&&&&&&&\\
                           &\ar@{.}@/_9pt/[dd] &&&&&&&&&&&&&&&&\\
                           & c_{m+1} \ar[ddd] &&&&&&&&&&&&&&&&\\
                           &&&  &&&&&&&&&&&&&&&& \\
                           &&&&& &&&&&&&&&&&&&& \\
                           & z \ar@{-}@/^/[dl] &&&&&&&&&&&&&&&\\
                           &&&&&&&&&&&&&&&&&&&}$$
\end{minipage}
\end{minipage}

\bigskip

If we have more than two cycles then  iterating the process we will have every cycle sharing an arrow with the  root cycle.

\end{obs}

It remains to see the particular case where the ray $\mathcal{R}$ does not have a linear part, that is, when the ray only has $m$-saturated cycles.

\begin{lema}\label{rel interna meter ciclos sin parte lineal}

The  gentle bound quiver

$$\SelectTips{eu}{10}\xymatrix@C=.4pc@R=.4pc{ &&&&&&\\
 & a \ar@{-}@/_/[ul] \ar[dd]&& && \\
 &\ar@{.}@/_9pt/[dd] & b_{m+1} \ar[dl]&& b_m \ar[ll] \ar@{.}@/^1.5pc/[dd]&&\\
 & b \ar[dd]  \ar[dr]&& \mathcal{C} &&\\
 && b_1 \ar[rr] && b_2 &&\\
 &c \ar@{-}@/^/[dl] &&&&\\
&&&&&&
}$$

is tilting-cotilting equivalent to a gentle bound quiver having the  cycle $\mathcal{C}$ sharing only one arrow with the root cycle.
\end{lema}

\begin{proof}
We have to consider the following cases:

\begin{enumerate}
 \item If there is no internal relation involving the vertex $a$,  apply the mutation $\sigma_b$ to get the bound quiver

 $$\SelectTips{eu}{10}\xymatrix@C=.4pc@R=.4pc{ && b_1 \ar[rr] && b_2 \ar@{.}@/^2pc/[dddd] &&\\
 & a \ar@{-}@/_/[ul] \ar[ur] && && \\
 &&&&&&\\
 & b\ar[uu] \ar[dd] & &&&&\\
 & \ar@{.}@/_1pc/[dd] & b_m \ar[lu]&& b_{m-1} \ar[ll]&&\\
 & b_{m+1} \ar[dd]&&&&&\\
 &&&&&&\\
 &c \ar@{-}@/^/[dl] &&&&\\
&&&&&&
}$$

 \item If there is no internal relation involving the vertex $c$,  apply the mutation  $\sigma'_b$ to obtain the bound quiver

$$\SelectTips{eu}{10}\xymatrix@C=.4pc@R=.4pc{ && &&&&\\
 & a \ar@{-}@/_/[ul] \ar[dd] && && \\
 &  \ar@{.}@/_1pc/[dd]&&&&&\\
 & b_1 \ar[dd] &&&\\
 & & b_2 \ar[rr] && b_3 \ar@{.}@/^2pc/[dddd] &&&&&\\
 & b \ar[ru] &&&&&\\
 &&&&&&\\
 & c \ar[uu] \ar@{-}@/^/[dl] & &&&&\\
&& b_{m+1} \ar[lu] && b_{m} \ar[ll]&&\\&&&&
}$$

 \item If there are relations involving the vertices $a$ and  $c$, to start, assume that at the vertex  $c$ there is no other  $m$-saturated cycle  $\mathcal{C'}$ attached. If this is the case then  we work with the  cycle $\mathcal{C'}$ instead of the cycle $\mathcal{C}$.\\
     Since the vertices  $a,b $ and  $c$ are on a non-oriented and a non-saturated cycle, we cannot have a path from  $c$ to $a$  with internal relations involving each vertex.  Thus, there is a vertex $c_n$ where the path starting at $c$ stops having internal relations. That is:

     $$\SelectTips{eu}{10}\xymatrix@C=.5pc@R=.7pc{c\ar[rr] & \ar@{.}@/^1pc/[rr] & c_1 \ar[rr] & & \cdots  \ar[rr] & \ar@{.}@/^1pc/[rr] & c_{n-2} \ar[rr] & \ar@{.}@/^1pc/[rr] & c_{n-1} \ar[rr] && c_n \ar@{-}[r] & }$$

 Now we move the cycle $\mathcal{C}$ to the vertex $c_n$ applying $n$ times the mutation $\sigma'_b$:

%

\begin{minipage}[c]{0.97\linewidth}
\begin{minipage}{0.43\linewidth}
$$\SelectTips{eu}{10}\xymatrix@C=.4pc@R=.4pc{
   &&&&&\\
 & a \ar@{-}@/_/[ul] \ar[dd]&& && \\
 &\ar@{.}@/_9pt/[dd] & b_{m+1} \ar[dl]&& b_m \ar[ll] \ar@{.}@/^1.5pc/[dd]&&\\
 & b \ar[dd]  \ar[dr]&& \mathcal{C} &&\\
 &\ar@{.}@/_9pt/[dd] & b_1 \ar[rr] && b_2 &&\\
 & c \ar[dd] &&&\\
 &\ar@{.}@/_9pt/[dd]&&&&&\\
 &c_1 \ar[dd] &&&\\
 &&&&&\\
 & \ar@{.}[dd]\\
 &&&&&\\
 &\ar[dd]&&&&\\
 &\ar@{.}@/_9pt/[dd]&&&&&\\
 & c_{n-1}\ar[dd] &&&&\\
 &&&&&\\
 & c_n \ar@{-}@/^/[dl] &&&&\\
&&&&&&}$$
\end{minipage}
\hfill
\begin{minipage}{0.1\linewidth}
$\substack{\stackrel{\sigma'_b}{\longmapsto}}$
\end{minipage}
\hfill
\begin{minipage}{0.43\linewidth}
$$\SelectTips{eu}{10}\xymatrix@C=.4pc@R=.4pc{
   &&&&&\\
 & a \ar@{-}@/_/[ul] \ar[dd]&& && \\
 &\ar@{.}@/_9pt/[dd]&&&&&\\
 &b_1 \ar[dd] &&&\\
 &\ar@{.}@/_9pt/[dd] & c \ar[dl]&& b_{m+1} \ar[ll] \ar@{.}@/^1.5pc/[dd]&&\\
 & b \ar[dd]  \ar[dr]&& \mathcal{C} &&\\
 &\ar@{.}@/_9pt/[dd] & b_2 \ar[rr] && b_3 &&\\
 &c_1 \ar[dd] &&&\\
 &&&&&\\
 & \ar@{.}[dd]\\
 &&&&&\\
 &\ar[dd]&&&&\\
 &\ar@{.}@/_9pt/[dd]&&&&&\\
 & c_{n-1}\ar[dd] &&&&\\
 &&&&&\\
 & c_n \ar@{-}@/^/[dl] &&&&\\
&&&&&&}$$
\end{minipage}
\end{minipage}

\bigskip

 and we are in a situation similar to that of case $2$.
\end{enumerate}
\end{proof}

\subsection{Moving the $m$-saturated cycles}
In order to get an algebra having the normal form we need to move closer together the $m$-saturated cycles that we have moved into the root cycle.


\begin{lema}\label{juntar todos los ciclos saturados}
Let $(Q,I)$ be a bound solar quiver without rays where each $m$-saturated cycle shares just one arrow with the root cycle. Then $(Q,I)$ is tilting-cotilting equivalent to a bound quiver such  that whenever two neighbouring $m$-saturated  cycles are connected by a path they end up connected by a vertex.
\end{lema}

\begin{proof}
Assume that   $\mathcal{C}$ and $\mathcal{B}$ are two neighbouring  $m$-saturated cycles such that the length of the path between them is at least one. If this path does not have a relation involving the first arrow  the quiver looks like one of the  quivers on the  left which, applying the sequence of mutations described below, changes to one of the right:

 \begin{center}
\begin{tabular}{lll}
  $\SelectTips{eu}{10}\xymatrix@C=.03pc@R=.9pc{ & c_m  \ar@{.}@/^1pc/[r] \ar[ddl]& c_{1} &&&  && && b_{m+1} \ar@{.}@/^1pc/[r] \ar@{-}[ddl]& b_2&&&&\\
   &  &&&&& &&&  & &&&\\
   e \ar[rrr] \ar@{-}[d]   &&& c \ar[uul]  \ar[rrr] & && d \ar@{.}[rr] && b \ar@{-}[rrr] &&& b_1\ar@{-}[uul] \ar@{-}[d] \\
   &&&&&&&&&&&
   }$
    &$\SelectTips{eu}{10}\xymatrix@C=.03pc@R=1pc{ & \\
\stackrel{\sigma'_d(\sigma'_{c_{m}}\sigma'_c)\cdots (\sigma'_{c_1}\sigma'_c)}{\longmapsto}\\
&\\
&\\
&\\
&}$ &

$\SelectTips{eu}{10}\xymatrix@C=.03pc@R=.9pc{ &&& & c_{m-1}  \ar@{.}@/^1pc/[r] \ar[ddl]& c_{1}   && && b_{m+1} \ar@{.}@/^1pc/[r] \ar@{-}[ddl]& b_2 &&&&\\
  &&&  &&&&&&&&&&\\
  c_m\ar[rrr] \ar@{-}[d] &&& c \ar[rrr]    &&&  d \ar[luu]\ar@{.}[rr] && b \ar@{-}[rrr] &&& b_1\ar@{-}[uul] \ar@{-}[d] \\
  &&& &&&&&&&&&&&
   }$

\end{tabular}
\end{center}

 \begin{center}
\begin{tabular}{lll}
  $\xymatrix@C=.06pc@R=.9pc{ & c_2 \ar@{.}@/^1pc/[r] & c_{m+1} \ar[rdd] &&&  && && b_{1} \ar@{.}@/^1pc/[r] & b_m \ar@{-}[ddr]&&&&\\
   &&&&&&&&&&&&&\\
   c_1  \ar@{-}[d] \ar[ruu]   &&& c \ar[lll]  \ar[rrr] & && d \ar@{.}[rr] && b \ar@{-}[ruu]  &&& b_{m+1} \ar@{-}[lll] \ar@{-}[d] \\
   &&&&&&&&&&&
   }$
    &$\SelectTips{eu}{10}\xymatrix@C=.6pc@R=1pc{ & \\
\stackrel{\sigma'_c}{\longmapsto}\\
&\\
&\\
&\\
&}$ &

$\xymatrix@C=.06pc@R=.9pc{ &&& & c_{2}  \ar@{.}@/^1pc/[r] & c_{m+1} \ar[rdd]  && && b_{1} \ar@{.}@/^1pc/[r] & b_m \ar@{-}[ddr] &&&&\\
  &&&  &&&&&&&&&&\\
  c_1\ar[rrr] \ar@{-}[d] &&& c \ar[ruu]    &&&  d \ar@{.}[rr]  \ar[lll] && b \ar@{-}[uur]&&& b_{m+1} \ar@{-}[lll] \ar@{-}[d] \\
  &&& &&&&&&&&&&&
   }$

\end{tabular}
\end{center}
Note that the cycles that are drawn without orientation can be oriented arbitrarily.
If instead,   the path has  a relation involving the first arrow  the quiver looks like one of the  quivers on the  left which, applying the sequence of mutations described below, changes to one of the right:

 \begin{center}
\begin{tabular}{lll}
  $\xymatrix@C=.06pc@R=.9pc{ & c_2 \ar@{.}@/^1pc/[r] & c_{m+1} \ar[ddr] & &&& &&& &&&  &b_1 \ar@{.}@/^1pc/[r] \ar@{-}[ddl] & b_m  \ar@{-}[ddr]& \\
  &&& &&& &&& &&& &&&&\\
  c_1 \ar@{-}[d] \ar[ruu] &&&  c \ar[lll] \ar[rrr] &&\ar@{.}@/^/[rr]& a \ar[rrr] &&& d \ar@{.}[rrr] &&& b \ar@{-}[rrr] &&& b_{m+1}\ar@{-}[d] \\
  &&&&&&&&&&&&&&&&}$

    &$\SelectTips{eu}{10}\xymatrix@C=.6pc@R=1pc{ & \\
\stackrel{\sigma'^2_c}{\longmapsto}\\
&\\
&\\
&\\
&}$ &

$\xymatrix@C=.06pc@R=.9pc{ &b_m \ar@{-}[ddl] \ar@{.}@/^1pc/[r] & b_1 \ar@{-}[ddr] &  &&& & a \ar[ddl]\ar@{.}@/^1pc/[r] & c_3 &  &&& &&& \\
&&& &&& &&& &&& &&&\\
b_{m+1}\ar@{-}[d]  &&& d \ar@{.}[rrr] \ar@{-}[lll] &&& b \ar[rrr]  &&& c \ar[uul] && \ar@{.}@/^/[rr] & c_2 \ar[lll] &&& c_1 \ar[lll]\ar@{-}[d]\\
&&&&&&&&&&&&&&&}$

\end{tabular}
\end{center}

\begin{center}
\begin{tabular}{lll}

$\SelectTips{eu}{10}\xymatrix@C=.06pc@R=.9pc{ & c_m \ar@{.}@/^1pc/[r]& c_1 & &&& &&& &&& & a_1\ar@{.}@/^1pc/[r] & a_m \ar[ddr] & \\
&&& &&& &&& &&& &&&\\
c_{m+1} \ar@{-}[d] \ar@{-}[uur]  \ar@{-}[rrr] &&& c \ar@{-}[uul] \ar@{.}[rrr] &&& d \ar[rrr] && \ar@{.}@/^/[rr] & a \ar[rrr] &&& b  \ar[ruu]&&& a_{m+1} \ar[lll]\ar@{-}[d] \\
&&& &&& &&& &&& &&&}$

    &$\SelectTips{eu}{10}\xymatrix@C=.6pc@R=1pc{ & \\
\stackrel{\sigma^2_b}{\longmapsto}\\
&\\
&\\
&\\
&}$ &

$\SelectTips{eu}{10}\xymatrix@C=.06pc@R=.9pc{ & c_m\ar@{-}[ddl] \ar@{.}@/^1pc/[r] & c_1 & &&& & a\ar@{.}@/^1pc/[r] & a_{m-1} \ar[ddr] &&&&&&\\
&&& &&& &&& &&& &&&\\
c_{m+1}\ar@{-}[uur] \ar@{-}[rrr]  \ar@{-}[d] &&& c\ar@{.}[rrr] \ar@{-}[uul] &&& d \ar[uur] &&& b \ar[rrr] \ar[lll]&& \ar@{.}@/^/[rr]& a_m \ar[rrr] &&& a_{m+1} \ar@{-}[d] \\
&&& &&& &&& &&& &&&}$

\end{tabular}
\end{center}

Thus, the length of the path between $\mathcal{C}$ and $\mathcal{B}$ decreases by 1 in the first case and by 2 in the second. If needed we can move the relations over the path using:

\begin{center}
\begin{tabular}{lll}
$\SelectTips{eu}{10}\xymatrix@C=.5pc@R=.5pc{ a_{j-2} \ar[rr] && a_{j-1} \ar[rr] & \ar@{.}@/^1pc/[rr]& a_j \ar[rr] && a_{j+1} \ar[rr] &&     }$
&
$\SelectTips{eu}{10}\xymatrix@C=.8pc@R=.9pc{
\stackrel{\sigma_{a_{j-1}}}{\longmapsto}}
$
&
$\SelectTips{eu}{10}\xymatrix@C=.5pc@R=.5pc{
 a_{j-2} \ar[rr] & \ar@{.}@/^1pc/[rr]& a_{j} \ar[rr] & & a_{j-1} \ar[rr] && a_{j+1} \ar[rr] &&    }$

\end{tabular}
\end{center}

and iterate the process.
\end{proof}

\subsection{Rays without union relations}
There remains the case where a ray is attached to the root cycle without an union relation. That is, through an $m$-saturated cycle.
An iteration of the following  lemma (and its dual) and  remark \ref{rel interna meter ciclos} allows to move such a ray into the root cycle.

\begin{lema}\label{rayo sin rel de union}

The  gentle bound quiver

 $$\SelectTips{eu}{10}\xymatrix@C=.4pc@R=.5pc{ && b_1 \ar[rr] && b_2 \ar@{.}@/^2pc/[dddd] &&\\
 & a \ar@{-}@/_/[ul] \ar[ur] && && \\
 &&&&&&\\
 &c\ar[uu] \ar@{-}@/^/[dl] &&&&\\
 &  & b \ar[lu] \ar[dd] && b_{m-1} \ar[ll]&&\\
 & &&&&&\\
 && b_m &&&&
 }$$


is tilting-cotilting equivalent to the gentle bound quiver

 $$\SelectTips{eu}{10}\xymatrix@C=.4pc@R=.5pc{ && b_1 \ar[rr] && b_2 \ar@{.}@/^2pc/[dddd] &&\\
 & a \ar@{-}@/_/[ul] \ar[ur] && && \\
 &&&&&&\\
 & b\ar[uu]  & &&&&\\
 &  & b_m \ar[lu]&& b_{m-1} \ar[ll]&&\\
 &c \ar@{-}@/^/[dl] \ar[uu] &&&&\\
&&&&&&
}$$

\end{lema}

\begin{proof}

 We apply $\sigma'_{b}$ to move the arrow $b\rightarrow b_m$ inside the root cycle.
\end{proof}

Observe that the proof still works if the vertices  $b_1,\cdots, b_i$ (where $1\leq i\leq m-1$) belong to the root cycle.


We are interested in reducing any solar algebra to an algebra having a  normal form. In these algebras every $m$-saturated cycle shares with the root cycle just one arrow. In light of the preceding result, we know that when we move onto the root cycle a ray having $m$-saturated cycles we do it in such a way that that condition holds. There remains to see the particular case where there are $m$-saturated cycles attached directly to the root cycle but sharing with it more than one arrow.


After remark \ref{ciclo q comparte mas de una flecha entonces hay flechas - 1 relaciones en el sentido contrario} we know that if an  $m$-saturated cycle shares with the root cycle $k-1$ arrows counterclockwise oriented  and  $m+2-k+1$ arrows clockwise oriented,  then we are in one of the following cases:

\begin{itemize}
  \item [(a)] there are at least $k-2$ clockwise internal relations or at least  $m+2-k$ counterclockwise internal relations.
  \item [(b)] there is another  $m$-saturated cycle sharing with the root cycle $k-1$ arrows clockwise oriented  and  $m+2-k+1$ arrows counterclockwise oriented.
\end{itemize}

If the condition (a) holds, then after moving some arrows and relations,  we can assume that we have the following situation:

$$\SelectTips{eu}{10}\xymatrix@C=.6pc@R=.8pc{ b_{k-1} \ar[rr] & \ar@{.}@/^/[rr] & b_{k-2} \ar[rr] &\ar@{.}@/^/[rr] & b_{k-3} \ar@{.}[rr]  & \ar@{.}@/^/[rr] & b_2 \ar[rr] & \ar@{.}@/^/[rr] & b_1 \ar[rr] && c_1 \ar[drrrr] &  & c_2 \ar[ll] &  & c_3 \ar[ll]  && c_{k-2}  \ar@{.}[ll] &  & c_{k-1} \ar[ll] && c_k \ar[ll] \\
 &&  &&   &&  && &&  &&&& c_{m+2} \ar@{.}[rr] && c_{k+1} \ar[urrrr]
}$$

where the vertices $b_{k-1}$ and  $c_k$ belong to the root cycle. Then, an iteration of the following lemma gives the claim.

\begin{lema}\label{cambiar sentido de un ciclo}



A gentle bound quiver as above is tilting-cotilting equivalent to the gentle bound quiver

$$\SelectTips{eu}{10}\xymatrix@C=.6pc@R=.8pc{ b_{k-1} \ar[rr] & \ar@{.}@/^/[rr] & b_{k-2} \ar[rr] & \ar@{.}@/^/[rr] & b_{k-3} \ar@{.}[rr]  &\ar@{.}@/^/[rr] & b_2 \ar[rr] &  & c_2 \ar[rr] && c_1 \ar[drr]  &  & c_3 \ar[ll]  && c_{k-2}  \ar@{.}[ll] &  & c_{k-1} \ar[ll] && c_k \ar[ll] \\
 &&  &&   &&  && &&  &&  b_1 \ar[rr] && c_{m+2} \ar@{.}[rr] && c_{k+1} \ar[urr]}$$


\end{lema}

\begin{proof}

We apply the sequence of mutations $\sigma_{c_2} \sigma_{c_1}$.
\end{proof}


If instead  $(b)$ holds, we can assume that both  $m$-saturated cycles share a vertex.  Then, an iteration of lemma  \ref{dos ciclos saturados se pueden pegar en el vertice que queramos} gives us the claim.

\section{Proof of the Main Theorem}

We start  providing the procedure to reduce any  $\tilde{\mathbb{A}}$-branched algebra to the normal form, defined in section \ref{normal forms}, using the local mutations from section \ref{mutaciones}. 

\begin{prop}\label{algoritmo}
 Let  $A$ be an  $\tilde{\mathbb{A}}$-branched algebra. Then, $A$ is tilting-cotilting equivalent to an oriented or non-oriented normal form.
\end{prop}

\begin{proof}
We can assume that $A$ is a solar algebra, then the result follows upon executing the following: \newline

{\bf Algorithm. } \begin{description}
 \item[Step 1] For each external union relation move the involved ray onto the root cycle removing the relation at the same time using  lemmas \ref{la flecha en el otro sentido}, \ref{relacion externa sin ciclos saturados}, and \ref{si tiene ciclos el rayo}.
 \item[Step 2] For each internal union relation move the involved ray into the root cycle using  lemmas \ref{rel interna meter flecha en sentido contrario} and  \ref{rel interna meter rayo} and  remark \ref{rel interna meter ciclos}.
\item[Step 3] For each ray without union relations  move the ray into the root cycle using  lemmas \ref{rel interna meter ciclos sin parte lineal} and \ref{rayo sin rel de union}.
 \item[Step 4] Move the  $m$-saturated cycles closer together  using lemma \ref{juntar todos los ciclos saturados}. 
 \item[Step 5] Move together the clockwise  $m$-saturated cycles:

\begin{center}
\begin{tabular}{lll}
$\xymatrix@C=.1pc@R=.9pc{ & a_m  \ar@{.}@/^1pc/[r] \ar[ddl]& a_{1} && b_1 \ar@{.}@/^1pc/[r] & b_m \ar[ddr] && c_m \ar@{.}@/^1pc/[r] \ar[ddl]& c_1&&&&\\
   &&&&&&&&&&&&&\\
   a\ar[rrr] \ar@{-}[d]   &&& b \ar[uul] \ar[uur] &&& c\ar[lll] \ar[rrr] &&& d\ar[uul] \ar@{-}[d] \\
   &&&&&&&&&&&
   }$
&
$\SelectTips{eu}{10}\xymatrix@C=.6pc@R=1pc{ & \\
\stackrel{\sigma'_c}{\longmapsto}\\
&\\
&\\
&\\
&}$
&
$\xymatrix@C=.1pc@R=.9pc{ & a_m  \ar@{.}@/^1pc/[r] \ar[ddl]& a_{1} && c_m \ar@{.}@/^1pc/[r] & c_1 \ar[ddr] && b_1 \ar@{.}@/^1pc/[r] \ar[ddl]& b_m &&&&\\
   &&&&&&&&&&&&&\\
   a\ar[rrr] \ar@{-}[d]   &&& b \ar[rrr]\ar[uul] \ar[uur] &&& c  &&& d\ar[uul] \ar[lll] \ar@{-}[d] \\
   &&&&&&&&&&&
   }$
\end{tabular}
\end{center}
\vspace*{-1cm}
Doing this we also move the counterclockwise cycles together.

 \item[Step 6] Move the clockwise arrows together:

\begin{center}
\begin{tabular}{lll}
$\SelectTips{eu}{10}\xymatrix@C=.8pc@R=.9pc{ a \ar[rr] && b && c \ar[rr] \ar[ll] && d
   }$
&
$\SelectTips{eu}{10}\xymatrix@C=.8pc@R=.9pc{
\stackrel{\sigma'_b}{\longmapsto}}
$
&
$\SelectTips{eu}{10}\xymatrix@C=.8pc@R=.9pc{
 a && b \ar[rr]  \ar[ll] && c \ar[rr] && d   }$

\end{tabular}
\end{center}

 \item[Step 7] Move the relations:

\begin{center}
\begin{tabular}{lll}
$\SelectTips{eu}{10}\xymatrix@C=.5pc@R=.5pc{ a_{j-2} \ar[rr] && a_{j-1} \ar[rr] & \ar@{.}@/^1pc/[rr]& a_j \ar[rr] && a_{j+1} \ar[rr] &&     }$
&
$\SelectTips{eu}{10}\xymatrix@C=.8pc@R=.9pc{
\stackrel{\sigma_{a_{j-1}}}{\longmapsto}}
$
&
$\SelectTips{eu}{10}\xymatrix@C=.5pc@R=.5pc{
 a_{j-2} \ar[rr] & \ar@{.}@/^1pc/[rr]& a_{j} \ar[rr] & & a_{j-1} \ar[rr] && a_{j+1} \ar[rr] &&    }$

\end{tabular}
\end{center}

 \item[Step 8] After iterating as many times as necessary steps 6 and 7, we can assume that the linear part (the one without $m$-saturated cycles) looks like:

     $$\SelectTips{eu}{10}\xymatrix@C=.4pc@R=.5pc{
   && c \ar[rrdd]  \ar[lldd] && \\
   & \ar@{.}@/^/[ddl] &&\ar@{.}@/_/[ddr]&\\
   . \ar[dd] &&&& . \ar[dd]\\
   \ar@{.}@/^/[dd] &&&& \ar@{.}@/_/[dd]\\
   . \ar[dd] &&&& . \ar[dd]\\
   \ar@{.}@/^/[dd] &&&& \ar@{.}@/_/[dd]\\
   .\ar@{.}[dd] &&&& . \ar@{.}[dd]\\
   &&&&\\
   &&&&
   }$$

 and remove the internal relations by pairs  (one in the clockwise direction with another in the counterclockwise direction):

\begin{center}
\begin{tabular}{lcccccl}
$\SelectTips{eu}{10}\xymatrix@C=.6pc@R=.6pc{
   && c \ar[rrdd]  \ar[lldd] && \\
   & \ar@{.}@/^/[ddl] &&\ar@{.}@/_/[ddr]&\\
   b_2 \ar[dd] &&&& b_1 \ar[dd]\\
   &&&&\\
   a_2 \ar@{.}[dd] &&&& a_1 \ar@{.}[dd]\\
   &&&&\\
   &&&&
   }$
&&&
$\SelectTips{eu}{10}\xymatrix@C=.6pc@R=1pc{ & \\
&\\
&\\
\stackrel{\tiny{\sigma'_{b_2}\sigma'_{b_1}\sigma'_c}}{\longmapsto}\\
&\\
&\\
&}$
&&&
$\SelectTips{eu}{10}\xymatrix@C=.6pc@R=.6pc{&& c \ar[rrdd]  \ar[lldd] && \\
   &  &&&\\
   b_1 \ar[dd] &&&& b_2 \ar[dd]\\
   &&&&\\
   a_2 \ar@{.}[dd] &&&& a_1 \ar@{.}[dd]\\
   &&&&\\
   &&&& }$

\end{tabular}
\end{center}

 Note that we do not need to have the same number of clockwise and counterclockwise relations. Thus, if $\alpha_h$ is the number of clockwise internal relations and  $\alpha_a$ the number of counterclockwise internal relations we will finally have $|\alpha_h-\alpha_a|$  relations in the direction of the bigger number.
\end{description}
\end{proof}

%
%
%
%
%

The normal forms are not necessarily   $m$-cluster tilted algebra of   type  $\tilde{\mathbb{A}}$. For instance,  the normal forms can have more than $m-1$ consecutive relations. The following corollary shows how to get an $m$-cluster tilted algebra of   type  $\tilde{\mathbb{A}}$ tilting-cotilting equivalent to the normal forms.

\begin{coro}\label{atildebranche-->der eq a m-cluster}
 Let $A$ be an  $\tilde{\mathbb{A}}$-branched algebra, then $A$ is  tilting-cotilting equivalent to an $m$-cluster tilted algebra of   type  $\tilde{\mathbb{A}}$.
\end{coro}

\begin{proof}
We know that  $A$ is tilting-cotilting equivalent to an algebra $N$ with the   normal form (oriented or not). Let $r$  be  the number of consecutive relations in  $N$. We can write $r=\a(m-1)+\b$ with $\b< m-1$. The condition $(c)$ from the definition of $\tilde{\mathbb{A}}$-branched algebra says that there are $r+1+\varepsilon$ arrows not belonging to any  $m$-saturated cycle, where

    $$\varepsilon =\left\{
              \begin{array}{ll}
                \alpha - 1, & \hbox{\text{if }  $\beta=0$ ;} \\
                \alpha, & \hbox{\text{if}  $\beta\neq 0$ .}
              \end{array}
            \right.$$

Then, using the procedure described in step 7 of the previous algorithm we can move the  relations in order to form $\a$ sets with  $m-1$ consecutive relations and  a final set with $\b$ consecutive relations. Moreover, we can separate each group by only  a vertex with no relations involved.  These algebras are always $m$-cluster tilted of  type  $\tilde{\mathbb{A}}$.
\end{proof}

One of the main results of this paper says that the converse also holds true, namely that the algebras derived equivalent to, the connected component of, $m$-cluster tilted algebras of type
$\tilde{\mathbb{A}}$ having the root cycle are precisely the $\tilde{\mathbb{A}}$-branched algebras.


 \begin{prop}\label{vale el reciproco}
  Let $B$ be a connected algebra with a root cycle derived equivalent to an   $m$-cluster tilted algebra of type $\tilde{\mathbb{A}}$. Then    $B$ is   $\tilde{\mathbb{A}}$-branched.
 \end{prop}

\begin{proof}

Since  $B$ is derived equivalent to an  $m$-cluster tilted algebra of  type $\tilde{\mathbb{A}}$, it follows that $B$ is gentle and derived equivalent to an  algebra $N=\widetilde{N}_{n_1,k_1,n_2,k_2,r}$ or $B_{k,n,t}$ where $|r|\equiv 0$ modulo $m$ for the first case and $n-1\equiv 0$ modulo $m$ (but $n-1\neq 0$) for the second case.

Assume that the  difference between the numbers of clockwise  and counterclockwise internal  relations for $B$ is $r'$. Then,  the algorithm of proposition  \ref{algoritmo} implies that  $B$ is derived equivalent to an  algebra $\widetilde{N}_{n_1',k_1',n_2',k_2',r'}$; if the root cycle is not oriented or an algebra  $B_{k',r'+1,t'}$; if it is oriented. Hence, if $B$ is derived equivalent to an algebra  $\widetilde{N}_{n_1',k_1',n_2',k_2',r'}$,  corollary  \ref{si dos formas normales son der eq entonces r y r' son cong} gives that $|r'|\equiv |r| $ modulo $m$ and since $N$ is derived equivalent to an   $m$-cluster tilted algebra of type $\tilde{\mathbb{A}}$, we conclude that $|r'|\equiv |r|\equiv 0 $ modulo $m$. Now, assume that   $B$ is derived equivalent to $B_{k',r'+1,t'}$. It remains to show that the root cycle of $B$  has at least one relation. This follows directly from the equality $\phi_B=\phi_{B_{k,n,t}}$.

Finally, write $r'=\alpha (m-1) + \beta$ (with $\beta< m-1$) with
    $$\varepsilon =\left\{
              \begin{array}{ll}
                \alpha - 1, & \hbox{\text{if }  $\beta=0$ ;} \\
                \alpha, & \hbox{\text{if}  $\beta\neq 0$ .}
              \end{array}
            \right.$$

            and assume that  $r'> 0$ (equivalently,  the relations are in the clockwise sense).  If we do not have at least  $r'+1+\varepsilon$ arrows on the root cycle clockwise oriented and  not belonging to any  $m$-saturated cycle,  then the $r'$  relations cannot be separated  into  sets with at most $m-1$ consecutive relations. In consequence, $B$ cannot be derived equivalent to an $m$-cluster tilted algebra of type $\tilde{\mathbb{A}}$.
\end{proof}
%
%
%
%
%
%
%
%
%

We are now able to state and prove the main theorem of this paper. Observe that Theorem A  in the Introduction is precisely  the equivalence between conditions $(a)$ and $(d)$.

\begin{teo}\label{teo principal}
Let  $A$ be a connected algebra with a root cycle. Then the following conditions are equivalent.
\begin{itemize}
  \item [(a)] $A$ is an  $\tilde{\mathbb{A}}$-branched algebra.
  \item [(b)] $A$ is tilting-cotilting equivalent to an algebra $\widetilde{N}_{n_1,k_1,n_2,k_2,r}$ or $B_{k,n,t}$.
  \item [(c)] $A$ is derived equivalent to an algebra $\widetilde{N}_{n_1,k_1,n_2,k_2,r}$ or $B_{k,n,t}$.
  \item [(d)] $A$ is derived equivalent to an  $m$-cluster tilted algebra of  type  $\tilde{\mathbb{A}}$.
\end{itemize}
\end{teo}

\begin{proof} \
\begin{enumerate}
 \item[\ ] (a) implies (b). This is proposition \ref{algoritmo}.
 \item[] (b) implies (c). This is immediate, since Happel showed in \cite{H88} that if two algebras are tilting-cotilting equivalent then  they are derived equivalent.
 \item[] (c) implies (d).  This follows  from corollary \ref{atildebranche-->der eq a m-cluster} since the normal forms are themselves  $\tilde{\mathbb{A}}$-branched algebras.
  \item[] (d) implies (a). Follows from proposition \ref{vale el reciproco}.
\end{enumerate}
\end{proof}

\section{Further consequences}

\subsection{Cluster tilted algebras of type $\widetilde{\mathbb{A}}$}
 Propositions \ref{phi de la forma normal orientada} and \ref{phi de la forma normal A tilde} together with Theorem \ref{teo principal} enable us to recover the classification of the algebras derived equivalent to cluster tilted algebras of type $\widetilde{\mathbb{A}}$, of Bobi\'nski and Buan  \cite[Theorem B]{BB10}.  Note however that we additionally obtained an explicit description of the bound quiver of such an algebra. The conditions characterizing those bound quivers are really easy to verify, for no computations are required, as in \cite[Theorem B]{BB10}, which we can easily recover.

\begin{coro}
Let $A$ be a gentle algebra. Then $A$ is derived equivalent to a cluster tilted algebra of type $\widetilde{\mathbb{A}}$ if and only if there exist natural numbers $m_1,m_2,p,q$ such that $p + m_1 > 0 $,  $q + m_2 > 0$  and

$$\phi_A=(m_1+m_2).(0,3)^*+(p+m_1,p)^*+(q+m_2,q)^*$$
\end{coro}

\begin{proof}
Assume $A$ is derived equivalent to a cluster tilted algebra of type $\widetilde{\mathbb{A}}$. From Theorem \ref{teo principal}, $A$ is derived equivalent to an oriented or non-oriented normal form, and $m=1$. Then, use Proposition \ref{phi de la forma normal A tilde} (or \ref{phi de la forma normal orientada}) with $m_1=k_1+r$, $m_2=k_2-r$, $p=n_1-k_1$ and $q=n_2-k_2$ (or $m_1=n-1$, $m_2=k+1-n$, $p=0$ and $q=n+t$, respectively).

On the other hand, assume $A=\K Q/I$ is gentle such that $\phi_A=(m_1+m_2).(0,3)^*+(p+m_1,p)^*+(q+m_2,q)^*$. Since $\phi_A(0,3)=m_1+m_2$, in $(Q,I)$ there are exactly $m_1+m_2$ $1$-saturated cycles.  In addition, \cite[Lemma 7.2]{BB10} gives that $A$ is an algebra with root, and since $\widetilde{\mathbb{A}}$-branched algebras for $m=1$ are exactly algebras with root we are done.
\end{proof}

\subsection{Derived equivalence classification}

Here we provide a complete classification of $m$-cluster tilted algebras of type $\widetilde{\mathbb{A}}$ up to derived equivalence and specializing to the case $m=1$  we recover the main result of Bastian \cite[Theorem 5.5]{Bastian09}.

Let $(Q,I)$ be the quiver of an $m$-cluster tilted algebra of type $\widetilde{\mathbb{A}}$. We define five parameters $s_1,s_2,k_1,k_2$ and $r$ for $(Q,I)$ as follows:

\begin{defi}
Let $s_1$ be the number of arrows which are not part of any $m$-saturated cycle and which fulfill one of the following conditions:

\begin{itemize}
  \item  [a)] These arrows are part of the root cycle and they are oriented in the counterclockwise direction.
  \item  [b)] These arrows belong to a ray attached to the root cycle by a counterclockwise internal union relation and this relation does not involve the arrows.
  \item  [c)] These arrows belong to a ray attached to the root cycle by a clockwise internal union relation and this relation  involve the arrows.
  \item  [d)] These arrows belong to a ray attached to the root cycle by a counterclockwise external union relation or a ray without union relations.
\end{itemize}

Let $k_1$ be the number of  $m$-saturated cycles  which fulfill one of the following conditions:

\begin{itemize}
  \item  [a)] These cycles share one arrow $\alpha$ with the root cycle and $\alpha$ is oriented in the counterclockwise direction.
  \item  [b)] These cycles belong to a ray attached to the root cycle by a clockwise internal union relation.
  \item  [c)] These cycles belong to a ray attached to the root cycle by a counterclockwise external union relation or a ray without union relations..
\end{itemize}

Similary we define the parameters $s_2$ and $k_2$   permuting the words  'counterclockwise' and 'clockwise'.\\

Let $r$ be the number of  clockwise internal relations less the number of counterclockwise internal relations. (Note that as in the Normal form  the relations are clockwise oriented if $r>0$ and are counterclockwise oriented instead).

\end{defi}

\begin{lema}\label{equivalencias entre Q}
Let $N=\widetilde{N}_{n_1,k_1,n_2,k_2,r}$ be a normal form. Assume that $r=mx+r'$, then $N$ is derived equivalent to $\widetilde{N}_{n_1,k_1+x,n_2,k_2-x,r'}$.
\end{lema}

\begin{proof}
Use Lemma \ref{cambiar sentido de un ciclo} to change the direction of  $x$ clockwise $m$-saturated cycles.
\end{proof}

Now, we are ready to show the main result about derived equivalence classification of cluster tilted algebras of type $\widetilde{\mathbb{A}}$.

\begin{teo}
Let $A$ and $B$ be two $m$-cluster tilted algebras of type $\widetilde{\mathbb{A}}$. Denote by  $s_1,s_2,k_1,k_2$ and $r$ ( or $s'_1,s'_2,k'_1,k'_2$ and $r'$ ) the parameters for $A$ (for $B$, respectively). Then, $A$ and $B$ are derived equivalent if and only if $s_1=s'_1$, $s_2=s'_2$, $k_1+k_2=k'_1+k'_2$ and $r-r'=m(k'_1-k_1)=m(k_2-k'_2) $.
\end{teo}

\begin{proof}
Let $N_A$ be  the normal form derived equivalent to $A$ and $N_B$  the normal form derived equivalent to $B$. Assume  $A$ is derived equivalent to $B$, then $N_A$ is derived equivalent to $N_B$ and hence  $\phi_{N_A}=\phi_{N_B}$.  Then using proposition \ref{phi de la forma normal A tilde} or \ref{phi de la forma normal orientada} we get the four equalities listed above.

Conversly, assume that the   parameters for  $A$ and $B$ satisfy the equations  $s_1=s'_1$, $s_2=s'_2$, $k_1+k_2=k'_1+k'_2$ and $r-r'=m(k'_1-k_1)=m(k_2-k'_2) $.  It is enough to show that $N_A$ and $N_B$ are derived equivalent. We can assume without loss of generality that the normal forms are non-oriented. Then, we know that $N_A=\widetilde{N}_{s_1+k_1, k_1, s_2+k_2, k_2, r}$.   Since $r-r'=m(k'_1-k_1)=m(k_2-k'_2)$, lemma \ref{equivalencias entre Q} implies that $N_A$ is derived equivalent to $\widetilde{N}_{s_1+k_1+x,k_1+x,s_2+k_2-x,k_2-x,r'}=\widetilde{N}_{s'_1+k'_1, k'_1, s'_2+k'_2, k'_2, r'}=N_B$.
\end{proof}

\subsection{Avella-Alaminos and Geiss derived invariant}

It was not known in general  whether or not $\phi$ is a complete invariant for gentle algebras. It was shown in \cite{Avella-Alaminos2008} that if $A$ and $A'$ two algebras whose quivers have Euler characteristic $1$ are such that $\phi_A = \phi_{A'}$, then $A$ and $A'$ are derived equivalent, and a similar result in \cite{AA07}. Further, Bobi\'nski and Buan showed in \cite{BB10} that $\phi$ is a complete invariant for cluster tilted algebras and Bustamante and Gubitosi in \cite{GubBust} proved that this is also the case for $m$-branched algebras.  Now we show that for the family of  $\tilde{\mathbb{A}}$-branched algebras, $\phi$ is not a complete invariant.

\begin{ejem}
 Fix $m=2$ and let $A$ and  $B$ be two   $\tilde{\mathbb{A}}$-branched algebras with the normal forms   $\widetilde{N}_{6,4,5,3,1}$ and $\widetilde{N}_{7,5,4,2,-1}$ respectively. Remember that, as usual, dotted lines mean relations, and  each oriented $4$-cycle is $2$-saturated.

\end{ejem}

\begin{figure}[H]
\begin{center}
\includegraphics[scale=.6]{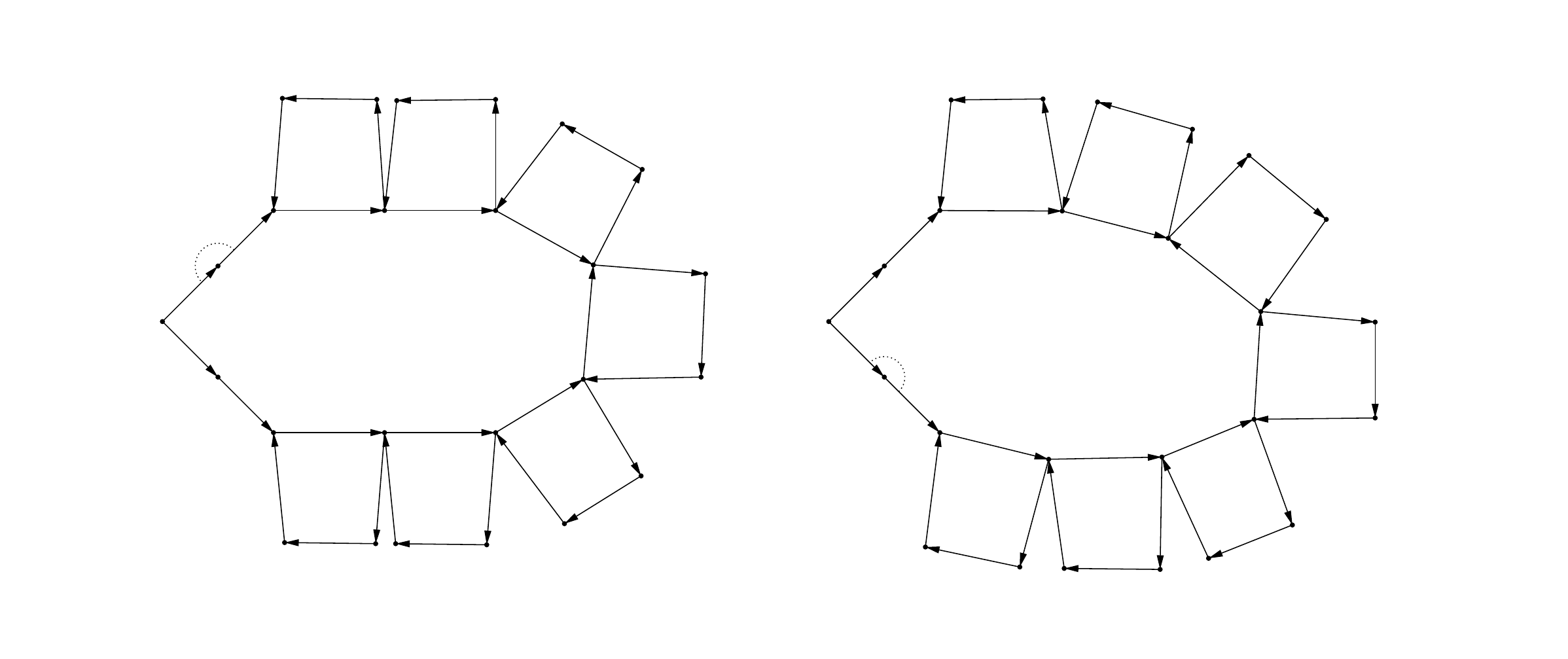}
\end{center}
\end{figure}

\vspace*{-0.80cm}
$$\widetilde{N}_{6,4,5,3,1} \hspace*{6cm} \widetilde{N}_{7,5,4,2,-1}$$

\medskip
Then $\phi_A=\phi_B$ but clearly $A$ and  $B$ are not derived equivalent.

\section*{Acknowledgements}
The author gratefully thanks Ibrahim Assem for several interesting and helpful discussions. This paper is part of the  author's Ph.D. thesis, done under the direction of Ibrahim Assem, she gratefully acknowledges financial support from the faculty of Sciences of the Universit\'e de Sherbrooke. 

\bibliographystyle{acm}

\bibliography{library}

\end{document}